\documentclass[11pt,openbib]{article}

\usepackage[T1]{fontenc}
\usepackage{amsmath}
\usepackage{amsfonts}
\usepackage{amssymb}
\usepackage{graphicx}

\newcommand{\Hh}{\mathbb{H}}

\newcommand{\R}{\mathbb{R}} 
\newcommand{\C}{\mathbb{C}}

\newcommand{\gl}{\mathop{\mathrm{gl}}\nolimits}
\newcommand{\Gl}{\mathop{\mathrm{Gl}}\nolimits}

\newcommand{\Oo}{\mbox{$\mathrm{O}$}}
\newcommand{\oo}{\mbox{$\mathrm{o}$}}
\newcommand{\spp}{\mathop{\mathrm{sp}}\nolimits}
\newcommand{\Spp}{\mathop{\mathrm{Sp}}\nolimits}

\newcommand{\setrule}{\, \rule[-4pt]{.5pt}{13pt}\, }

\newcommand{\onehalf}{\mbox{$\frac{\scriptstyle 1}{\scriptstyle 2}\,$}} 
\newcommand{\spann}{\mathop{\rm span}\nolimits} 

\newcommand{\ttfrac}[2]{\mbox{$\frac{{\scriptstyle #1}}{{\scriptstyle #2}}$}}

\newcommand{\vvee}{\mbox{\tiny $\vee $}}

\allowdisplaybreaks

\begin{document}
 
\hspace{.25in}{\Large \parbox[t]{4.5in}{\textbf{Adjoint orbits of real affine Lie algebras} } }\\ 
\par \noindent \hspace{2.5in}Richard Cushman  
\vspace{.5in}
\addtocounter{footnote}{1}
\footnotetext{printed: \today} 

\noindent {\Large \textbf{Introduction} \medskip
%%%%%%%%%%%%%%%%

\normalsize In this paper we find a representative of each orbit of the adjoint 
action of a real affine classical group of its Lie algebra. These orbits are 
not determined by the usual Jordan invariants of eigenvalues and 
block sizes, but require a noneigenvalue parameter called a modulus. 
This modulus was discovered by Weyl \cite{weyl} for the 
affine general linear group. The cases for the real affine orthogonal 
and symplectic groups have been treated in \cite{cushman-vanderkallen} and 
\cite{cushman}, respectively. The arguments in \cite{burgoyne-cushman} and 
\cite{cushman-vanderkallen} have been the motivation for the approach 
used here. Many essential details are new. We treat each real classical 
affine group as an isotropy group of a complex classical group with an 
anti-linear automorphism of order $2$. This ensures 
that we do not have to work with vector spaces defined over the quaternions. 
\medskip 

Here is a list of the contents of each section. \S 1 defines a real classical 
group. \S 2 defines the notion of type as in \cite{burgoyne-cushman} 
and distinguished type as in \cite{cushman-vanderkallen}. Our main result 
is: every distinguished type is the sum of an indecomposable nilpotent distinguished type and a sum of indecomposable types, which is unique up to reordering of 
the sumands, see table 2 on page \ref{pg-table2}. \S 3 classifies indecomposable nilpotent distinguished types. 
\S 4 lists the isotropy groups of the real classical groups that are isomorphic 
to a real affine classical group. 

\section{The real classical groups}
%%%%%%%%%%%%%%%%%%%%%%%%

Let $\widetilde{V}$ be a complex vector space of dimension $n$ and 
let  $\Gl ( \widetilde{V} )$ be the group of invertible complex linear mappings 
of $\widetilde{V}$ into itself. Let $\widetilde{\tau }$ be a nondegenerate 
symmetric or alternating complex bilinear form on $\widetilde{V}$. Set 
$G(\widetilde{V}, \widetilde{\tau}) = \{ g \in \Gl (\widetilde{V}) \setrule \, 
g^{\ast}\tau = \tau \} $.
Recall that $(g^{\ast }\tau )(\widetilde{u}, \widetilde{v}) 
= \tau \big( g(\widetilde{u}), g(\widetilde{v}) \big) $ for every $\widetilde{u}$, 
$\widetilde{v} \in \widetilde{V}$. Then $G( \widetilde{V}, \tau ) = 
\Oo (\widetilde{V}, \widetilde{\tau } )$ if $\widetilde{\tau }$ is symmetric 
and $\Spp (\widetilde{V}, \widetilde{\tau } )$ is $\widetilde{\tau }$ is 
alternating. $G$ is a \emph{complex classical group} if it is equal to 
$\Gl (\widetilde{V})$, $\Oo (\widetilde{V} , \widetilde{\tau })$ or 
$\Spp (\widetilde{V}, \widetilde{\tau })$. \medskip 

We now define a real classical group. Suppose that $\widetilde{\sigma }: 
\widetilde{V} \rightarrow \widetilde{V}$ is an anti-linear map, that is, 
$\widetilde{\sigma }(\widetilde{u} + \widetilde{v}) = \widetilde{\sigma }(\widetilde{u}) 
+ \widetilde{\sigma }(\widetilde{v})$ for every $\widetilde{u}$, $\widetilde{v} 
\in \widetilde{V}$ and for every $\alpha \in \C$ one has 
$\widetilde{\sigma }(\alpha \widetilde{v}) = 
\overline{\alpha }\, \widetilde{\sigma }(\widetilde{v})$ for every 
$\widetilde{v} \in \widetilde{V}$. Moreover for every $g \in G$ suppose that \medskip  

\noindent either
\par \parbox[t]{4in}{${\widetilde{\sigma }}^{\ast }_{\pm}(g) = 
{\widetilde{\sigma }}^{-1}_{\pm }g {\widetilde{\sigma }}_{\pm}$ for every 
$g \in G$, where $(\widetilde{\sigma })^2_{\pm } = 
\pm \, {\mathrm{id}}_{\widetilde{V}}$ 
and ${\widetilde{\sigma }}^{\ast}_{\pm}\tau = \overline{\tau }$;}
\par \noindent or
\par \parbox[t]{4in}{$\widetilde{\sigma}: V \rightarrow V^{\ast }$ such that 
${\widetilde{\sigma }}^{\ast }(g) = {\widetilde{\sigma }}^{-1}g^T\widetilde{\sigma }$, 
where $g^T:V^{\ast } \rightarrow V^{\ast }$ is defined by 
$\big( g^T ({\widetilde{u}}^{\ast }) \big) (\widetilde{v}) = 
{\widetilde{u}}^{\ast }\big( g(\widetilde{v}) \big) $ for every 
${\widetilde{u}}^{\ast } \in {\widetilde{V}}^{\, \ast }$ 
and every $\widetilde{v} 
\in \widetilde{V}$. Here $g \in \Gl (\widetilde{V} )$. Let 
${\widetilde{\tau }}_{\ast }(\widetilde{u}, \widetilde{v}) = 
\big( \widetilde{\sigma }(\widetilde{u})\big) (\widetilde{v})$. Then 
${\widetilde{\tau}}_{\ast }$ is a hermitian form on $\widetilde{V}$.} \medskip 

\noindent Then $\widetilde{\sigma }$ is an automorphism of $G$ of order $2$. Let 
$G(\widetilde{V}, \widetilde{\sigma }, \widetilde{\tau}) = 
\{ g \in G \setrule \, {\widetilde{\sigma }}^{\ast }(g) = g \, \, \& \, \, 
g^{\ast}\tau = \tau \} $. Then $G(\widetilde{V}, \widetilde{\sigma }, \widetilde{\tau })$ 
is a \emph{real classical group}. For an explicit description of the anti-linear mappings, see appendix $1$. \bigskip 

\par \noindent \hspace{.25in}\begin{tabular}{cll}
\multicolumn{1}{c}{Complex group} & \multicolumn{1}{c}{Real group} & 
\multicolumn{1}{c}{Standard notation} \\ \hline
\rule{0pt}{12pt}$\Gl (\widetilde{V})$ & 
$\Gl (\widetilde{V}, {\widetilde{\sigma }}_{+})$ & \rule{10pt}{0pt}$\Gl (n, \R ) $ \\
 & $\Gl (\widetilde{V}, {\widetilde{\sigma }}_{-})$,\, \mbox{$n$ even} & 
\rule{10pt}{0pt}${\mathrm{U}}^{\ast }(n)$ \\
& $\Gl (\widetilde{V}, {\widetilde{\tau}}^{\, (p)}_{\ast } )$, \, $0 \le p \le 
\big[ \frac{n}{2} \big] $ & \rule{10pt}{0pt}$\mathrm{U}(n-p,p)$ \\
$\Oo (\widetilde{V}, \widetilde{\tau })$ & $\Oo (\widetilde{V}, \widetilde{\tau }, 
{\widetilde{\sigma}}^{\, (p)}_{+})$, \, $0 \le p \le \big[ \frac{n}{2} \big]$ & 
\rule{10pt}{0pt}$\Oo (n-p,p)$ \\
& $\Oo (\widetilde{V}, \widetilde{\tau }, {\widetilde{\sigma }}_{-})$, \, \mbox{$n$ even} & \rule{10pt}{0pt}${\Oo }^{\ast }(n)$ \\
$\Spp (\widetilde{V}, \widetilde{\tau })$ & $\Spp (\widetilde{V}, \widetilde{\tau }, 
{\widetilde{\sigma }}_{+})$, \, \mbox{$n$ even} & \rule{10pt}{0pt}$\Spp (n, \R )$ \\
& $\Spp (\widetilde{V}, \widetilde{\tau}, {\widetilde{\sigma}}^{\, (p)}_{-}) $, \, 
$0 \le p \le \big[ \frac{n}{4} \big] $ & \rule{10pt}{0pt}$\Spp (n-p,p)$
\end{tabular} \bigskip 

\par \noindent \hspace{.25in}\parbox[t]{4.5in}{Table 1. The real classical groups. See appendix $1$ for the definition of ${\tau }^{(p)}_{\ast }$ and  
${\widetilde{\sigma }}^{\, (p)}_{\pm}$. } \bigskip 

\noindent A real classical group $G(\widetilde{V}, \widetilde{\sigma }, \widetilde{\tau })$ 
is a Lie group with Lie algebra $g(\widetilde{V}, \widetilde{\sigma }, 
\widetilde{\tau})$ $ = \{ X \in \gl (\widetilde{V}) \setrule \, 
{\widetilde{\sigma }}^{\ast }X = X \, \, \& \, \, \widetilde{\tau }(X\widetilde{u}, 
\widetilde{v}) + \widetilde{\tau }(\widetilde{u}, \widetilde{v}) =0 \, \, 
\mbox{for every $\widetilde{u}$, $\widetilde{v} \in \widetilde{V}$} \} $. 

\section{Classification of adjoint orbits}
%%%%%%%%%%%%%%%%%%%

Let $G(\widetilde{V}, \widetilde{\sigma}, \widetilde{\tau })$ be a 
real classical group. Suppose that $v^0$ is a nonzero vector in 
$\widetilde{V}$, which lies in an eigenspace of the complex linear 
map $\overline{\widetilde{\sigma}}: \widetilde{V} \rightarrow \widetilde{V}$ 
defined by $\overline{\widetilde{\sigma }}(\alpha \widetilde{v}) = 
\alpha \widetilde{\sigma }(\widetilde{v})$ for every $\alpha \in \C$ and 
every $\widetilde{v} \in \widetilde{V}$. Note that 
${\overline{\widetilde{\sigma }}}^{\, 2}_{\pm} = \pm \, {\mathrm{id}}_{\widetilde{V}}$. 
In addition, suppose that $v^0$ is $\widetilde{\tau }$ isotropic, that is, 
$\tau (v^0,v^0) =0$. Let 
${G(\widetilde{V}, \widetilde{\sigma }, \widetilde{\tau })}_{v^0} 
= \{ g \in G(\widetilde{V}, \widetilde{\sigma }, \widetilde{\tau }) \setrule \, 
g(v^0) = v^0 \}$. The isotropy group ${G(\widetilde{V}, \widetilde{\sigma }, \widetilde{\tau })}_{v^0}$ is a closed subgroup of the Lie group 
$G(\widetilde{V}, \widetilde{\sigma }, \widetilde{\tau })$ and hence is a Lie 
group. Its Lie algebra is 
${g(\widetilde{V}, \widetilde{\sigma }, \widetilde{\tau })}_{v^0} = 
\{ X \in g(\widetilde{V}, \widetilde{\sigma }, \widetilde{\tau }) \setrule \, 
Xv^0 = 0 \}$. The adjoint action of ${G(\widetilde{V}, \widetilde{\sigma }, \widetilde{\tau })}_{v^0}$ 
on its Lie algebra ${g(\widetilde{V}, \widetilde{\sigma }, \widetilde{\tau })}_{v^0}$ 
is given by 
\begin{displaymath}
{G(\widetilde{V}, \widetilde{\sigma }, \widetilde{\tau })}_{v^0} \times 
{g(\widetilde{V}, \widetilde{\sigma }, \widetilde{\tau })}_{v^0} \rightarrow 
{g(\widetilde{V}, \widetilde{\sigma }, \widetilde{\tau })}_{v^0} :  
(P,X) \longmapsto PXP^{-1}. 
\end{displaymath}

We begin the classification of adjoint orbits  by defining 
the concept of an indecomposable type and an indecomposable 
distinguished type. Let $W$ be a $\tau $-nondegenerate $\sigma$ invariant 
subspace of the finite dimensional 
complex vector space $V$. If $Y \in g(W, \sigma |W ; \tau |W)$, the Lie algebra 
of $G(W, \sigma |W; \tau |W)$, then $(Y, W; \sigma |W, {\tau }|W)$ 
is a \emph{pair}. Two pairs $(Y, W; \sigma |W, $ ${\tau }|W)$ and 
$(Y', W'; $ $ \sigma |W', {\tau }|W')$ are \emph{equivalent} if there is a bijective 
complex linear mapping $P: W \rightarrow W'$ such that 
$Y' = PYP^{-1}$, ${\sigma }' = P\sigma P^{-1}$, and $P^{\ast }{\tau }' = 
\tau $. Being equivalent is an equivalence relation on the set of pairs. An 
equivalence class of pairs is a \emph{type}, denoted by $\Delta $, 
which is represented by the pair $(Y, W; \sigma |W, {\tau }|W)$. We define the 
\emph{dimension} of $\Delta $, denoted $\dim \Delta $, by $\dim W$ and 
the \emph{index} of $\Delta $, denoted $\mathrm{ind}\, \Delta $, by 
the number of negative eigenvalues of the Gram matrix 
$\big( \tau (v_i,v_j) \big) $ of $\tau $ with respect to the basis 
$\{ v_1, \ldots , v_{\dim W} \}$ of $W$. Let $Y = S+N$ be the Jordan 
decomposition of $Y$ into a semisimple linear map $S$ and a commuting 
nilpotent map $N$, both of which lie in $g(W, \sigma , \tau )$. To see this 
we argue as follows. Since $\sigma N {\sigma }^{-1}$ and 
$\sigma S {\sigma }^{-1}$ are commuting nilpotent and semisimple linear 
maps, whose sum is $\sigma Y {\sigma }^{-1} = Y$ and the Jordan 
decomposition of $Y$ is unique, it follows that $N = \sigma N {\sigma }^{-1}$ 
and $S = \sigma S {\sigma }^{-1}$. So $N$ and $S$ commute with 
$\sigma $. A similar argument shows that $N^{\ast }\tau = \tau $ and 
$S^{\ast }\tau = \tau $. Thus $S$ and $N$ lie in $g(W, \sigma , \tau )$. 
There is a unique nonnegative integer $h$ such that $N^h \ne 0$ but 
$N^{h+1}= 0$. We call $h$ the 
\emph{height} of of the type $\Delta $ and denote it by $\mathrm{ht}(\Delta )$. 
$\mathrm{ht}(\Delta )$ does not depend on the choice of representative 
of $\Delta $. We say that the type $\Delta $, represented by the 
pair $(W, \sigma , \tau )$, is \emph{uniform} if $NW = \ker N^hW$. Suppose that 
$W = W_1 \oplus W_2$, where $W_i$ are proper, $Y$ invariant subspaces of 
$V$, which are $\sigma $ invariant, $\tau $ nondegenerate, and $\tau $ 
orthogonal. Then $\Delta $ is the \emph{sum} of two types ${\Delta }_i$, 
represented by $(Y|W_i, W_i; \sigma |W_i, \tau |W_i)$. We write 
$\Delta = {\Delta }_1+ {\Delta }_2$. The type $\Delta $ is \emph{indecomposable} 
if it cannot be written as the sum of two types. In other words, 
the pair $(Y,W; \sigma , \tau )$ represents an indecomposable type $\Delta $ if 
there is no proper $Y$ and $\sigma $ invariant subspace of $W$ on which 
$\tau $ is nondegenerate. From 
\cite[prop.3 p.343]{burgoyne-cushman} it follows that an indecomposable type 
is uniform. The classification of indecomposable types is given in 
\cite{burgoyne-cushman}. \medskip 

Next we introduce the concept of a triple. Let $\overline{\sigma}$ be the complex 
linear mapping of $V$ into itself associated to the anti-linear mapping 
$\sigma $ defined by setting $\overline{\sigma }(\alpha v) = 
\alpha \sigma (v)$ for every $v \in V$ and every $\alpha \in \C$. Again let 
$W$ be a $\sigma $ (and thus $\overline{\sigma }$) invariant subspace of 
$V$, which is $\tau $ nondegenerate. We say that 
$(Y,W, v^0; \sigma , \tau )$ is a \emph{triple} if 1) $v^0$ is a nonzero vector in $W$. 
If $\sigma $ is involved, then $v^0$ is an eigenvector of $\overline{\sigma }$ corresponding to the eigenvalue 
$\lambda $, that is $\pm 1$ if $\sigma  = {\sigma }_{+}$ or $\pm \mathrm{i}$ if 
$\sigma = {\sigma }_{-}$; 2) if $\tau $ or ${\tau }_{\ast}$ is involved, then 
$v^0$ is $\tau $ isotropic, that is, $\tau (v^0,v^0) = 0$ or ${\tau }_{\ast }(v^0,v^0) =0$ ; 3) The linear map $Y \in g(W; \sigma , \tau )$ and $Yv^0 =0$.  If $v$ is a nonzero vector in $V$ which satisfies condition 1) and 2) of the definition of a triple, then we call $v$ \emph{special}. Because $S$ and $N$ are polynomials in $Y$ with complex coefficients and $Yv^0 =0$, it follows that $Sv^0 = Nv^0 = 0$. So $S$, $N \in {g(W, \sigma , \tau )}_{v^0}$.  We say that two triples $(Y,W, v^0; \sigma , \tau )$ and $(Y',W', (v')^0; {\sigma }', {\tau }')$ are \emph{equivalent} if there is a 
complex bijective linear mapping $P:W \rightarrow W'$ such that 
$Y' = PYP^{-1}$, $Pv^0 = (v')^0$, ${\sigma }' = P\sigma P^{-1}$, and 
$P^{\ast }{\tau }' =\tau $. Note that the vector $(v')^0$ is a special vector in 
$V'$, since 
\begin{align*}
\overline{{\sigma }'}\big( (v')^0 \big) & = {\sigma }'\big( (v')^0 \big) = 
P \sigma \big( P^{-1}(v')^0 \big) = P \big( \sigma (v^0) \big) \\
& = P\big( \overline{\sigma }(v^0) \big) = P(\lambda v^0), \, \, \mbox{where $\lambda = \pm 1$ or $\pm \mathrm{i}$} \\
& = \lambda P(v^0) = \lambda (v')^0,  
\end{align*}
and 
\begin{displaymath}
{\tau }'\big( (v')^0, (v')^0 \big) = {\tau }'(Pv^0, Pv^0) = 
(P^{\ast }{\tau }')(v^0,v^0) = \tau (v^0,v^0) = 0. 
\end{displaymath} 
Being equivalent is an equivalence relation on the 
set of triples. We call an equivalence class of triples a 
\emph{distinguished type}, which we denote by $\underline{\Delta }$. 
We say that $\underline{\Delta}$ is represented by the triple $(Y, W, v^0; \sigma , 
\tau )$. Suppose that $W = W_1 \oplus W_2$, where $W_i$ are proper, 
$\sigma $ and $Y$ invariant, $\tau $ nondegenerate, and $\tau $ orthogonal 
subspaces of $V$ with $v^0$ a $\tau $ isotropic special vector in $W_1$. Then 
$(W_1, Y|W_1, v^0; \sigma |W_1, \tau |W_1)$ is a triple, whose corresponding distinguished type is 
${\underline{\Delta }}_1$. Let the pair $(W_2, Y|W_2; \sigma |W_2, \tau |W_2)$ 
represent the type ${\Delta }_2$. We say that the distinguished 
type $\underline{\Delta }$, represented by the triple $(Y|W, W, 
v^0; {\sigma }|W, {\tau }|W)$, is the \emph{sum} of the distinguished type 
${\underline{\Delta }}_1$ and the type ${\Delta }_2$, which we write as 
$\underline{\Delta } = {\underline{\Delta }}_1 + {\Delta }_2$. If $\underline{\Delta }$ 
can not be written as the sum of a distinguished type and and type, then 
$\underline{\Delta }$ is an \emph{indecomposable} distinguished type. 
In other words, the triple $(Y, W, v^0; \sigma , \tau )$ represents an 
indecomposable distinguished type if there is no proper, $Y$ and $\sigma $ 
invariant, $\tau $ nondegenerate subspace of $W$, which contains $v^0$. \medskip 

The goal of the next few sections is to prove \medskip 

\noindent \textbf{Theorem 2.1} Every distinguished type is the sum 
of an indecomposable nilpotent distinguished type and a sum of 
indecomposable types, which is unique up to reordering of the summands. \medskip 

The proof of this theorem requires classifying indecomposable distinguished 
types. Recall that indecomposable types have been classified in 
\cite{burgoyne-cushman}. Theorem 3.1 solves the conjugacy class 
problem for the Lie algebra ${g(V, \sigma , \tau )}_{v^0}$ under the 
adjoint action of ${G(V, \sigma , \tau )}_{v^0}$. Indeed, there is a one to one correspondence between equivalence classes triples and orbits of the adjoint 
action of ${G(V, \sigma , \tau )}_{v^0}$ on its Lie algebra 
${g(V, \sigma , \tau )}_{v^0}$. \medskip 

Suppose that $\underline{\Delta }$ is a distinguished type, which is 
represented by the triple $(Y, W, v^0; \sigma , \tau )$. We say that 
$\underline{\Delta }$ has \emph{distinguished height} $h$ if there is 
a special vector $w\in W$ such that $Y^hw = v^0$. We denote the 
distinguished height of $\underline{\Delta }$ by $\mathrm{dht}(\underline{\Delta })$. 
\medskip 

\noindent \textbf{Lemma 2.2} Suppose that $\underline{\Delta } = 
{\underline{\Delta }}' + \Delta $. Then $\mathrm{dht}(\underline{\Delta }) = 
\mathrm{dht}({\underline{\Delta }}')$.  \medskip 

\noindent \textbf{Proof.} First we show that $\mathrm{dht}(\underline{\Delta }) = 
\mathrm{dht}({\underline{\Delta }}')$. Suppose that $(Y, W, v^0; \sigma , \tau )$ 
is a triple which represents $\underline{\Delta }$ and that $W = 
W_1 \oplus W_2$, where $W_i$ are proper, $Y$ and $\sigma $ invariant, 
$\tau $ nondegenerate, $\tau $ orthogonal subspaces of $W$ with 
$v^0 \in W_1$. The triple $(Y|W_1, W_1, v^0; \sigma , \tau )$ represents 
${\underline{\Delta }}'$ and the pair $(Y|W_2, W_2; \sigma , \tau )$ represents 
the type $\Delta $. Suppose that $\mathrm{dht}({\underline{\Delta }}') = h'$. 
Then there is a special vector $w' \in W_1$ such that $Y^{h'}w' = v^0$. 
Since $W_1 \subseteq W$ we get $\mathrm{dht}(\underline{\Delta }) \ge h'$. 
Because $\mathrm{dht}(\underline{\Delta }) = h$, there is a special vector 
$w \in W$ such that $Y^hw = v^0$. But $W = W_1 \oplus 
W_2$. So we may write $w = w_1 + w_2$, where $w_i \in W_i$ and $w_{\ell }$ is 
a special vector. Since $W_i$ is $Y$ invariant, one has $v^0 = Y^hw = Y^hw_1 + Y^hw_2$, where $Y^hw_i \in W_i$. By construction $v^0 \in W_1$. Therefore 
$Y^hw_1 = v^0$. Consequently, $h\le  \mathrm{dht}({\underline{\Delta }}') = 
h'$. So $h = h'$. Note that 
$\dim {\underline{\Delta }}' < \dim \underline{\Delta }$. \medskip 

Let $\underline{\Delta }$ be a distinguished type of distinguished height 
$h$, which is represented by the triple $(Y,W, v^0; \sigma , \tau )$. If 
$\tau $ is \emph{not} involved in the triple representing $\underline{\Delta }$, 
let $\mu (\underline{\Delta })$ be the set of all $\lambda \in {\C} \setminus \{ 0 \}$ 
such that $\lambda Y^hw = v^0$ for some special vector $w \in W$. 
If $\tau $ is involved, let $\mu (\underline{\Delta })$ be the set of all 
$\tau (w, v^0) \in \C$ such that $w$ is a special vector in $W$ with $Y^hw = v^0$. We call $\mu (\underline{\Delta })$ the \emph{set of parameters} of 
the distinguished type $\underline{\Delta }$.\medskip 

\noindent \textbf{Lemma 2.3} If $\underline{\Delta } = {\underline{\Delta }}' + 
\Delta $ is a distinguished type of distinguished height $h$, then 
$\mu (\underline{\Delta }) = \mu ({\underline{\Delta }}')$. \medskip 

\noindent \textbf{Proof.} Suppose that $\tau $ is not involved in the 
triple $(Y,W, v^0; \sigma , \tau )$ representing $\underline{\Delta }$, and 
that $\lambda \in \mu ({\underline{\Delta }}')$. Then there is a special vector 
$w \in W_1$ such that $\lambda Y^hw = v^0$, since the distinguished height 
of ${\underline{\Delta }}'$, represented by the triple $(Y_1, W_1, 
v^0; \sigma |W_1, \tau |W_1)$,  is $h$ by lemma 2.2. But $W_1 \subseteq W$. 
So $\lambda \in \mu (\underline{\Delta })$, by definition. Hence 
$\mu ({\underline{\Delta }}') \subseteq \mu (\underline{\Delta })$. Next we prove 
the reverse inclusion. If $\lambda \in \mu (\underline{\Delta })$, there is a 
special vector $w \in W$ such that $\lambda Y^hw = v^0$. But $\underline{\Delta } = 
{\underline{\Delta }}' + \Delta $. So there are $Y$ and $\sigma $ invariant 
subspaces $W_i$ which are $\tau $ orthogonal such that $W = W_1 
\oplus W_2$ with $v^0 \in W_1$. Write $w = w_1+w_2$. Then  
$w_i$ are special vectors in $W_i$ with $\lambda Y^hw = \lambda Y^hw_1 +\lambda Y^hw_2$. But $Y^hw_i \in W_i$, since $W_i$ are $Y$ invariant and 
$\lambda Y^hw= v^0$. So $\lambda Y^hw_1 = v_0$ since $v^0 \in W_1$ by hypothesis, that is, $\lambda \in \mu ( {\underline{\Delta }}' )$. Thus $\mu (\underline{\Delta }) \subseteq \mu ( {\underline{\Delta } }')$, which shows that 
$\mu ( {\underline{\Delta }}' ) = \mu (\underline{\Delta })$ as desired.
\par Now suppose that $\tau $ is involved in the 
triple $(Y,W, v^0; \sigma , \tau )$ representing $\underline{\Delta }$. Since 
$W_1 \subseteq W$ it follows from the definition of the set of parameters that 
$\mu ( {\underline{\Delta } }' ) \subseteq \mu ( \underline{\Delta} )$. 
Suppose that there is a special vector $w \in W$ with $Y^hw = v^0$ 
such that $\tau (w, v^0) \notin \mu ( {\underline{\Delta } }' )$. Write 
$w = w_1+w_2$ with $w_i \in W_i$, which are special vectors. Then $Y^h w_1 = v^0$. Since $\underline{\Delta } = { \underline{\Delta } }' + \Delta $ the subspace $W_1$ is $\tau $ orthogonal to the subspace $W_1$ and $v^0 \in W_1$. So 
\begin{displaymath}
\tau (w, v^0) = \tau (w_1, v^0) + \tau (w_2, v^0) = \tau (w_1, v^0) .
\end{displaymath}
But $\tau (w_1, v^0) \in \mu ({\underline{\Delta }}')$. This is a contradiction. 
Hence $\mu ({\underline{\Delta }}') = \mu (\underline{\Delta })$. 
\hfill $\square $ \medskip 

\noindent \textbf{Lemma 2.4} Suppose that $\underline{\Delta }$ is a distinguished 
type. Then we may write $\underline{\Delta }$ $ = 
{\underline{\Delta }}' + \Delta$, where the distinguished type 
${\underline{\Delta }}' $ is indecomposable and nilpotent. \medskip 
 
\vspace{-.15in}\noindent \textbf{Proof.} If ${\underline{\Delta }}$ is not indecomposable, we find another distinguished type ${\underline{\Delta }}''$ of the same distinguished height and set of parameters as $\underline{\Delta }$ and a type 
${\Delta }$, where $\dim \Delta >0$, such that 
${\underline{\Delta }} = {\underline{\Delta }}'' + {\Delta }$. Because 
$\dim {\underline{\Delta }} > \dim {\underline{\Delta }}''$ after a finite number of repetitions, we obtain a distinguished type 
${\underline{\Delta }}'$, which we can not write as the sum of a distinguished type and a type $\widetilde{\Delta }$, namely, $\underline{\Delta } = 
{\underline{\Delta }}' + \widetilde{\Delta }$. In other words, 
${\underline{\Delta }}'$ is an indecomposable distinguished type, 
which by lemmas 2.2 and 2.3 has the same distinguished height and set of 
parameters as $\underline{\Delta }$.  
\par We now show that the indecomposable distinguished type 
$\underline{\Delta }'$, represented by the triple 
$(Y|W, W, v^0; \sigma , \tau )$ is nilpotent. Let $W_0$ be 
the generalized eigenspace of $Y|W$ corresponding to the 
eigenvalue $0$. Then $v^0 \in W_0$, \linebreak 
because $Yv^0 =0$. Moreover, $W_0$ is $Y$ and $\sigma $ invariant 
and is $\tau $ nondegenerate. On 
$W_0$ the linear map $Y$ is nilpotent. Since $\underline{\Delta }'$ 
is indecomposable, it follows that the triple $(Y|W, W, v^0; \sigma , \tau )$ is equal to the triple $(Y|W_0, W_0, v^0; \sigma , \tau )$. Hence the indecomposable 
distinguished type $\underline{\Delta }'$  is nilpotent. 
\hfill $\square $ \medskip 

Lemma 2.4 proves theorem 2.1 except for the uniqueness assertion. To 
prove uniqueness we need the classification of indecomposable nilpotent 
distinguished types given in the next section. 

\section{Classification of nilpotent indecomposable \\
\rule{25pt}{0pt}distinguished types}
%%%%%%%%%%%%%%%%%

In this section we classify nilpotent indecomposable distinguished types. \medskip 

First we prove \medskip 

\noindent \textbf{Proposition 3.1} Let $\underline{\Delta }$ be an 
indecomposable nilpotent distinguished type of distinguished height 
$h$ and set of parameters $\mu (\underline{\Delta })$. Suppose that 
$\underline{\Delta }$ is represented by the triple $(Y, W, v^0; \sigma , \tau )$, 
where $Y \in g(W,\sigma ,\tau )$ is nilpotent and $Yv^0 =0$. Then 
$\underline{\Delta }$ is uniform and $\dim \ker Y = 1$ or $2$. \medskip 

\noindent \textbf{Proof.} Suppose that $\tau $ is not involved in the 
triple defining $\underline{\Delta }$. Since $\mathrm{dht}(\underline{\Delta }) = h$, 
there is a special vector $w \in W$ and $\lambda \in \mu (\underline{\Delta })$ 
such that $\lambda Y^hw = v^0$. Look at the subspace $\widetilde{W} = 
{\spann }_{\C}\{ w, Yw, \ldots , Y^hw= v^0 \}$. Since $Y^{h+1}|W =0$, the 
subspace $\widetilde{W}$ is $Y$ invariant and contains the vector $v^0$. 
$\widetilde{W}$ is also $\sigma $ invariant, 
because $w$, being a special vector, is an eigenvector of the 
complex linear map $\overline{\sigma }$, which commutes with $Y$. Hence 
$\widetilde{W}$ is $\overline{\sigma }$ invariant and hence $\sigma $ invariant. 
But $\underline{\Delta }$ is indecomposable. So $\widetilde{W} = W$. 
Hence $\underline{\Delta }$ is uniform because $Y$ has only one Jordan block 
on $W$. Thus $\dim \ker Y =1$. \medskip    

Now suppose that $\tau $ is involved in the triple defining $\underline{\Delta }$. 
There are two cases. \medskip 

\noindent \textsc{Case} $\mathbf{1}$. $\mu (\underline{\Delta }) \ne \{ 0 \}$. 
Since $\mathrm{dht}(\underline{\Delta }) = h$ and $\mu (\underline{\Delta }) \ne 
\{ 0 \} $, there is a special vector $w \in W$ such that $Y^hw = v^0$ and 
$\tau (w, v^0) = \mu \ne 0$. The subspace $\widetilde{W} = {\spann }_{\C}\{ w, Yw, 
\ldots , Y^hw = v^0 \}$ contains $v^0$ and is both $Y$ and $\sigma $ invariant. 
Let $T$ be the $(h+1) \times (h+1)$ matrix, whose $ij^{\mathrm{th}}$ entry 
is $\tau (Y^{i-1}w, Y^{j-1}w)$. $T$ has nonzero entries on its anti-diagonal 
since 
\begin{align*}
\tau (Y^{i-1}w, Y^{j-1}w) & = (-1)^{i-1}\tau (w, Y^{i+j-2}w), \, \, 
\mbox{since $Y \in g(W,\sigma ,\tau )$} \\
&\hspace{-.75in} = (-1)^{i-1}\tau (w, Y^hw) , \, \, \mbox{since $i +j = h+2$ on the anti-diagonal} \\
&\hspace{-.75in} = (-1)^{i-1}\tau (w, v^0) = (-1)^{i-1} \mu \ne 0. 
\end{align*}
$T$ has zero entries below the anti-diagonal because $Y^{\ell }w =0$ when 
$\ell > h$. Thus $\det T \ne 0$. In other words, $\tau $ is nondegenerate 
on $\widetilde{W}$. But $\underline{\Delta }$ is indecomposable. Hence 
$\widetilde{W} = W$. Consequently, $\underline{\Delta }$ is uniform, 
because $Y$ is one Jordan block. Hence $\dim \ker Y =1$.  \medskip 

\noindent  \textsc{Case} $\mathbf{2}$. $\mu (\underline{\Delta }) = \{ 0 \}$. 
By the results of \cite[p.343]{burgoyne-cushman} the type $\Delta $, 
represented by the pair $(Y, W; \sigma , \tau )$, which underlies the 
distinguished type $\underline{\Delta }$, is the sum ${\Delta }_1+ \cdots + 
{\Delta }_k$ of indecomposable types. Here ${\Delta }_{\ell }$ has height 
$h_{\ell }$ and is represented by the pair $(Y|W_{\ell } , W_{\ell }, \sigma |W_{\ell }, 
\tau |W_{\ell })$ for $1 \le \ell \le k$. Without loss of generality we may assume 
that $h_1 \le h_2 \le \cdots \le h_k \le h$. Then $W = \sum^k_{\ell } \oplus 
W_{\ell }$, where $W_{\ell }$ are $Y$ and $\sigma $ invariant and pairwise 
$\tau $ orthogonal with $\tau |W_{\ell }$ nondegenerate. Write 
$v^0 = \sum^k_{\ell =1} v^0_{\ell }$ with $v^0_{\ell } \in W_{\ell }$. Let 
$j$ be the smallest integer with $1 \le j \le k$ such that $v^0_j \ne 0$. 
Then $Yv^0_j =0$, because $0 = Yv^0 = \sum^k_{\ell =1}Yv^0_{\ell }$ and 
$Yv^0_{\ell } \in W_{\ell }$ for all $1 \le \ell \le k$. There is a vector $z \in W_j$ 
such that $\tau (z, v^0_j) \ne 0$, because $\tau |W_j $ is nondegenerate. The 
following argument shows that we may assume that $z$ is a special vector. 
Since the subspace $W_j$ is $\sigma $, and hence $\overline{\sigma }$, 
invariant, we have $W_j = W^1_j \oplus W^2_j$, where $W^1_j$ is the 
eigenspace of $\overline{\sigma }|W_j$ corresponding to the nonzero eigenvalue 
$\lambda $ and $W^2_j$ is its eigenspace corresponding to the eignevalue 
$-\lambda $. Recall that $\lambda = \pm 1$ or $\pm \mathrm{i}$, since 
$(\overline{\sigma }|W_j)^2 = \pm 1$. Write 
$z = z_1+z_2$, where $z_i \in W^i_j$. Suppose that $\tau (z_i, v^0_j) =0$ for 
$i =1$ and $2$. Then $\tau (z,v^0_j) = \tau (z_1,v^0_j) + \tau (z_2, v^0_j) = 0$, 
which contradicts the definition of the vector $z$. Suppose that 
$\tau (z_1, v^0_j) \ne 0$. Otherwise interchange $z_1$ and $z_2$ and 
rename $z_1$ as $z$. Since $W_j$ is $\tau $ orthogonal to 
$W_{\ell }$ for every $\ell \ne j$, it follows 
that $\tau (z, v^0) = \tau (z, v^0_j) \ne 0$. Moreover, $Y^hz \ne 0$, for 
suppose that $Y^nz =0$, then 
\begin{displaymath}
0=\tau (Y^hz,w) =(-1)^h\tau (z, Y^hw) = (-1)^h\tau (z, v^0) \ne 0.
\end{displaymath}
Suppose that $h_j < h$. Then $W' = \sum^k_{\ell = j+1}\oplus W_{\ell } \ne \{ 0 \}$. By hypothesis that $\underline{\Delta }$ has distinguished height $h$, there is 
a special vector $w' \in W'$ such that $Y^hw' = v^0$. Since $W'$ is 
$Y$ invariant, $Y^hw' \in W'$. So $v^0 \in W'$. This implies $v^0_j =0$, 
because $v^0 = v^0_j + v'$ for some $v' \in W'$. This contradicts the definition 
of the vector $v^0_j$. Hence $h_j =h$. Since $\underline{\Delta }$ has distinguished height $h$, there is a special vector $w \in W$ such that $Y^hw = v^0$. Consider 
the subspace
\begin{displaymath}
\widetilde{W} = {\spann }_{\C} \{ w, Yw, \ldots , Y^hw = v^0; \, z, Yz, \ldots , Y^hz \} .
\end{displaymath}
Clearly $\widetilde{W}$ is $Y$ invariant and contains the vector $v^0$. 
$\widetilde{W}$ is also $\sigma $ invariant, because the vectors $z$ and 
$w$ are special. From the definition of the vectors $z$ and $w$ it follows that 
$\tau (z, Y^hw) = \tau (z, v^0) \ne 0$ and $\tau (w, Y^hw) = \tau (w, v^0) =0$, 
since $\mu (\underline{\Delta }) = \{ 0 \} $. Look at the $2(h+1) \times 2(h+1)$ matrix 
$T = ${\tiny $\left( \begin{array}{c|c}
\tau (Y^{i-1}w, Y^{j-1}w) & \tau (Y^{i-1}w, Y^{j-1}z) \\ \hline 
\tau (Y^{i-1}z, Y^{j-1}w) & \tau (Y^{i-1}z, Y^{j-1}z ) \end{array} \right) $,}
where $i+j = h+2$ and $1 \le i, j \le h+1$. Then 
\begin{displaymath}
T= \mbox{\footnotesize $\left( \begin{array}{ccc|ccc}
\mbox{\large $\ast$} & & 0 & \mbox{\large $\ast $}  & & +  \\ 
& 0 & & & + & \\
0 & & \mbox{\Large $0$} & + & & \mbox{\Large $0$} \\ \hline 
\mbox{\large $\ast $} & & + & \mbox{\large $\ast $} & & \ast  \\
& +& & & \ast  & \\
+ & & \mbox{\Large $0$} & \ast & & \mbox{\Large $0$} \end{array} \right)$.}
\end{displaymath}
where the entries below the anti-diagonal of each block of $T$ vanish, 
because $\tau (Yu, v) = - \tau (u, Yv)$ for every $u$, $v \in \widetilde{W}$ 
and $Y^{\ell }|W = 0$ for every $\ell >h$. The entries on the anti-diagonal 
of the upper right hand and lower left hand block are nonvanishing 
because $\tau (Y^{i-1}z, Y^{j-1}w) = (-1)^{i-1}\tau (z, Y^{i+j-2}w) = 
(-1)^{i-1}\tau (z, Y^hw) \ne 0$ and $\tau (Y^{i-1}w, Y^{j-1}z) = 
(-1)^{j-1}\tau (Y^hw, z) \ne 0$. The entries on the anti-diagonal of the upper 
left hand block vanish, because $\tau (w, Y^hw) =0$. To 
compute the determinant of $T$ we expand 
by minors of the $(h+2)^{\mathrm{nd}}$ row and $(h+1)^{\mathrm{st}}$ column. 
Removing this row and column gives a matrix of the same form as $T$ with
one fewer row and column. Since the entries on the anti-diagonal of the lower 
left hand block of $T$ are nonzero, we find that the determinant of $T$ is 
the product of these entries and the determinant of the upper right hand 
block. Hence $\det T$ is nonzero. Thus $\tau |\widetilde{W}$ is 
nondegenerate. Since $\underline{\Delta }$ is indecomposable, it follows that 
$\widetilde{W} = W$. Consequently, $\underline{\Delta }$ is uniform, 
because $Y|W$ is the sum of two Jordan blocks of the same size, and 
$\dim \ker Y =2$. \hfill $\square$ \medskip 

By proposition 3.1 the nilpotent indecomposable distinguished type 
$\underline{\Delta }$ of distinguished height $h$ and set of parameters 
$\mu (\underline{\Delta })$ is uniform. Thus the reduced type 
$\overline{\overline{\underline{\Delta }}}$, represented by the tuple 
$(\overline{\overline{Y}}, \overline{\overline{W}} = W/YW, 
\overline{\overline{(v^0)}}; \overline{\overline{\sigma }}, \overline{\overline{\tau }})$ 
exists. One has $v^0 \in NW$, because $v^0 \in \ker Y$ and $NW = \ker Y^h|W$. 
Hence $\overline{\overline{(v^0)}} =0$. Since $NW = \ker Y^h|W$ and 
$Y^{h+1}|W =0$, the induced map 
\begin{displaymath}
\overline{\overline{Y}}: \overline{\overline{W}} \rightarrow \overline{\overline{W}}: 
\overline{\overline{w}} = w +YW \mapsto 
\overline{\overline{Y}}\, \overline{\overline{w}} = Yw + YW = YW, 
\end{displaymath}
vanishes identically. Thus the reduced type $\overline{\overline{\Delta}}$ is semisimple and is represented by the pair $(0,\overline{\overline{W}}; 
\overline{\overline{\sigma}}, \overline{\overline{\tau }})$.  Here 
$\overline{\overline{\sigma }}:\overline{\overline{W}} \rightarrow \overline{\overline{W}}$ is the anti-linear map 
induced by $\sigma $ with ${\overline{\overline{\sigma }}}^{\, 2} = \pm 1$ and 
$\overline{\overline{\tau }}: \overline{\overline{W}} \rightarrow 
\overline{\overline{W}}$ is the bilinear form 
$\overline{\overline{\tau }}(\overline{\overline{u}}, 
\overline{\overline{v}}) = \tau (u, Y^hv)$ for $u$, $v \in W$. The form 
$\overline{\overline{\tau }}$ is well defined, because for 
$\overline{\overline{u}} = u +NW$ and $\overline{\overline{v}} = v +NW$
\begin{align}
\overline{\overline{\tau}} (u + YW, v+YW) & = \tau (u+YW, Y^hv), \, \, 
\mbox{since $Y^{h+1}|W =0$}  \notag \\
& = (-1)^h \tau (Y^h(u + YW), v), \, \, \mbox{since $Y\in g(W, \sigma ,\tau )$} 
\notag \\
& = (-1)^h \tau (Y^hu, v) = \tau (u, Y^hv). \notag 
\end{align}
Also $\overline{\overline{\tau }}$ is nondegenerate, for if 
$0 = \overline{\overline{\tau }}(\overline{\overline{u}}, \overline{\overline{v}})$ 
for every $\overline{\overline{u}} \in \overline{\overline{W}}$, then 
$0 = \tau (u, Y^hv)$ for every $u \in W$. So $Y^hv =0$, since $\tau $ is nondegenerate on $W$. Hence $\overline{\overline{v}}=0$. Note that 
$\overline{\overline{\tau }}$ is symmetric if $h$ is even and $\tau $ is symmetric or $h$ is odd and $\tau $ is alternating. Otherwise $\overline{\overline{\tau}}$ is alternating. \medskip 

Let $\Delta $ be the uniform nilpotent type of height $h$, represented by the pair 
$(Y,W; \sigma , \tau )$, whose reduced type is the semisimple type 
$\overline{\overline{\Delta }}$, represented by the pair $(0,\overline{\overline{W}}; 
\overline{\overline{\sigma }}, \overline{\overline{\tau }})$. \medskip 

\noindent \textbf{Corollary 3.1A} The nilpotent type $\Delta $ associated to the indecomposable nilpotent distinguished type $\underline{\Delta }$ of distinguished height $h$ and set of parameters $\mu (\underline{\Delta })$ is 
indecomposable. \medskip 

\noindent \textbf{Proof.} By \cite[prop.5 p.343]{burgoyne-cushman} we can 
write $\overline{\overline{\Delta }} = {\overline{\overline{\Delta }}}_1 + \cdots +{\overline{\overline{\Delta }}}_k$, where ${\overline{\overline{\Delta }}}_{\ell }$ is 
an indecomposable semisimple type represented by the pair 
$(0,{\overline{\overline{W}}}_{\ell }; 
\overline{\overline{\sigma }}|W_{\ell}, \overline{\overline{\tau }}|W_{\ell })$ for 
$1 \le \ell \le k$, see appendix 3. Let ${\Delta }_{\ell }$ be the nilpotent type determined by the semisimple type 
${\overline{\overline{\Delta }}}_{\ell }$. ${\Delta }_{\ell }$ is 
represented by the pair $(Y|W_{\ell }, W_{\ell }; \sigma |W_{\ell }, \tau |W_{\ell})$, 
see appendix 2. By \cite[prop.3 p.343]{burgoyne-cushman} the type 
${\Delta }_{\ell }$ is uniquely determined by ${\overline{\overline{\Delta }}}_{\ell }$ and the height $h_{\ell }$ and is indecomposable. Hence $\Delta = {\Delta }_1 + \cdots + {\Delta }_k$, where $W = \sum^k_{\ell =1} \oplus W_{\ell }$. Since 
$\underline{\Delta }$ has distinguished 
height $h$, there is a special vector $w \in W$ such that $Y^hw =v^0$. 
Write $w= \sum^k_{\ell =1}w_{\ell }$ with $w_{\ell } \in W_{\ell}$ and 
$v^0 = \sum^k_{\ell =1}v^0_{\ell }$ with $v^0_{\ell } \in W_{\ell }$. Then 
$\sum^k_{\ell =1}v^0_{\ell } = v^0 = Y^hw = \sum^k_{\ell =1}Y^hw_{\ell }$. 
Since the subspace $W_{\ell }$ is $Y$ invariant, it follows that 
$v^0_{\ell } = Y^hw_{\ell }$ for all $1 \le \ell \le k$. Since $w$ is a special 
vector in $W$, there is a nonzero $\lambda \in \C$ such that the vector $w$ is 
an eigenvector of the complex linear map $\overline{\sigma }$ corresponding 
to the eigenvalue $\lambda $. So 
\begin{displaymath}
\sum^k_{\ell =1} \lambda w_{\ell } = \lambda w = \overline{\sigma }(w) = 
\sum^k_{\ell =1}\overline{\sigma }(w_{\ell}).
\end{displaymath}
But $\overline{\sigma }(w_{\ell }) \in W_{\ell }$ for all $1 \le \ell \le k$, 
since $W_{\ell }$ is a subspace of $W$ which is $\sigma $ invariant and 
hence $\overline{\sigma }$ invariant. In other words, $w_{\ell }$ is 
a special vector in $W_{\ell}$ for $1 \le \ell \le k$. Let $j$ be the smallest 
integer $1 \le j \le k$ such that $v^0_j \ne 0$. Let $W' = 
\sum^k_{\ell = j+1}\oplus W_{\ell }$. 
\par Suppose that $\tau $ is not involved in the triple defining the distinguished 
type $\underline{\Delta }$. If $j=k$, then $Y^hw_k =v^0_k =v^0$. Set $w'' =w_k$ 
and go to equation (\ref{eq-3zero}). Suppose that there is a nonzero vector 
$v' \in W'$ such that $v^0 = v^0_j + v'$. Let $w' = \sum^k_{\ell = j+1} w_{\ell }$. 
Then $w'$ is a special vector in $W'$ and 
\begin{displaymath}
Y^hw' = \sum^k_{\ell =1} Y^hw_{\ell } = \sum^k_{\ell =j+1}v^0_{\ell }
=\sum^k_{\ell =j}v^0_{\ell } -v^0_j  = v^0 - v^0_j = v'. 
\end{displaymath}
Let $w'' = w_j + w'$. Then $w''$ is a special vector in $W_j \oplus W'$ such 
that $Y^hw'' = Y^hw_j + Y^hw' = v^0_j +v' = v^0$. Let 
\begin{equation}
\widehat{W} = {\spann }_{\C}\{ w'', Yw'', \ldots , Y^hw'' =v^0 \}. 
\label{eq-3zero}
\end{equation}
Then $\widehat{W} \subseteq W_j \oplus W'$ is a $Y$ and $\sigma $ invariant subspace of $W$, which contains $v^0$. Because $\underline{\Delta }$ is indecomposable, $\widehat{W} = W$. But this contradicts the fact that $W_j \oplus W'$ is a proper subspace of $W$. Thus $v'=0$, which implies $W' = \{ 0 \}$ 
and thus $W = W_j$. So $\Delta = {\Delta }_j$, which is indecomposable. 
\par Suppose that $\tau $ is involved in the triple representing the distinguished 
type $\underline{\Delta }$. Suppose that $\tau (w'', Y^hw'') = \mu \ne 0$. The $(h+1) \times (h+1)$ matrix $\big( \tau (Y^iw'', Y^j w'') \big) $ of $\tau $ with respect 
to the basis $\{ w'', \, Yw'', $ $\ldots , Y^hw'' \} $ of $\widehat{W}$ is 
\begin{displaymath}
T = \mbox{\footnotesize $\left( \begin{array}{ccc} \mbox{\Large $\ast $}&  & \mu \\
&  -\mu &  \\  
& & \\
\mu & & \mbox{\Large $0$}   
\end{array} \right) $}, 
\end{displaymath}
since $\tau (Y^iw'', Y^jw'') = 0$ if $(i-1)+(j -1)>h$. But $\det T \ne 0$. So 
${\tau }|\widehat{W}$ is nondegenerate. Hence $\widehat{W} = W$, 
because the distinguished type $\underline{\Delta }$ is indecomposable. This contradicts the fact that $\widehat{W}$ is a proper subspace of $W$. Hence 
the vector $v'=0$. So $\Delta  = {\Delta }_j$, 
which is indecomposable. 
\par Suppose that $\tau (w'', Y^hw'') =0$. Then either $\tau $ is symmetric or $\tau = {\tau }_{\ast }$ and $h$ is odd, or $\tau $ is alternating and $h$ is even. So 
$\widehat{W}$ is a $\tau $ isotropic subspace of $W$. Here  
$\tau $ is nondegenerate bilinear form on $W$, because for $1 \le \ell \le k$ one 
has  $\tau |W_{\ell }$ is nondegenerate, since ${\Delta }_{\ell }$ is a 
type, and the subspaces $W_{\ell}$ are pairwise $\tau $ orthogonal. Hence 
there is a nonzero vector $z \in W_j \oplus W'$ such that 
$\tau (z, Y^hw'') =\mu \ne 0$ and $\tau (z,Y^hz) = \nu \ne 0$. We may 
suppose that the vector $z$ is special, for we may write 
$z = z_1 +z_2$, where $z_i \in (W_j \oplus W') \cap W^i$ and $W^i$ is 
an eigenspace of $\overline{\sigma }|W$. Now 
$0 \ne \tau (z, Y^hw) = \tau (z_1,Y^hw) + \tau (z_2, Y^hw)$ implies that 
one of the summands is not zero, say $\tau (z_1,Y^hw) \ne 0$. Otherwise 
interchange $z_1$ and $z_2$ and rename $z_1$ to be $z$. Then 
$z$ is a special vector. Consider the subspace 
\begin{displaymath}
W^{\vvee} = \spann \{ z, \, Yz, \ldots , Y^hz; \, w'', \, Yw'', \ldots , Y^hw'' = v^0 \}  
\subseteq W_j \oplus W'. 
\end{displaymath}
Then $W^{\vvee}$ is $Y$ invariant. $W^{\vvee}$ is $\sigma $ invariant because the 
vectors $w''$ and $z$ are special. Also we have 
\begin{displaymath}
\tau (Y^iz, Y^jz)  = \left\{ \mbox{\footnotesize $\begin{array}{cl}
(-1)^j\nu , & \mbox{if $i+j =h$} \\
0, &  \mbox{if $i+j >h$,}  
\end{array} $} \right.    
\tau (Y^iw'', Y^jz)  = \left\{ \mbox{\footnotesize $\begin{array}{cl}
(-1)^j \mu , & \mbox{if $i+j =h$} \\
0, &  \mbox{if $i+j >h$,} \end{array} $} \right.  
\end{displaymath}  
and $\tau (Y^iw'', Y^jw'')  = 0 $. 
Thus the matrix $T$ of $\tau $ on $W^{\vvee}$ has nonvanishing entries on 
the antidiagonal of each block, zero entries below the 
antidiagonal of each block, and a zero lower right hand block, that is, 
\begin{displaymath}
T = \mbox{\footnotesize $\left( \begin{array}{ccc|ccc}
\mbox{\large $\ast $} & & + & \mbox{\large $\ast $} & & + \\
& + &  &  & + &  \\ 
+ & & \mbox{\Large $0$} & + & &\mbox{\Large $0$} \\ \hline
\mbox{\large $\ast $} & & + & & & \\
& + & & & \mbox{\Large $0$} & \\
+ & &\mbox{\Large $0$} & & & \end{array} \right) $}
\end{displaymath}
Because all the entries of $T$ below the main anti-diagonal vanish and 
all the entries on the main anti-diagonal are nonzero the  
determinant of $T$ is nonvanishing. So ${\tau }|W^{\vvee}$ is nondegenerate. 
Hence $W^{\vvee} = W$, because $\underline{\Delta }$ is indecomposable. 
This contradicts the fact that $W^{\vvee}$ is a proper subspace of $W$. 
Thus the vector $v' =0$. So $W' = \{ 0 \}$, that is, $W = W_j$. Hence $\Delta  = {\Delta }_j$, which is indecomposable. \hfill $\square $ \medskip 

We now classify nilpotent indecomposable distinguished types 
$\underline{\Delta }$ of distinguished height $h$ and set of parameters 
$\mu (\underline{\Delta })$. Suppose that $\underline{\Delta }$ is represented by the triple $(Y,W,v^0; \sigma , \tau )$. From proposition 3.1 it follows that 
the dimension of the reduced semisimple type 
$\overline{\overline{\underline{\Delta }}}$ is 
$1$ or $2$.  The results of this classification are given in 
table 2. We argue case by case. \bigskip 

\noindent \hspace{-.25in}\begin{tabular}{ccccc}
\multicolumn{1}{c}{Lie algebra} & 
\multicolumn{1}{c}{$\overline{\overline{\underline{\Delta }}}$} & 
\multicolumn{1}{c}{$\underline{\Delta }$} & 
\multicolumn{1}{c}{$\dim \overline{\overline{W}}$} & 
\multicolumn{1}{c}{Conditions} \\ \hline
\rule{0pt}{11pt} $\gl (V, {\sigma }_{+})_{v^0}$ & ${\Delta }_0(0)$ & 
$\underline{{\Delta }_h(0), \, \lambda \in \R \setminus \{ 0 \} }$ & $1$ 
& A \\
\rule{0pt}{11pt} $\gl (V, {\sigma }_{-})_{v^0}$ & ${\Delta }_0(0,0)$ & 
$\underline{{\Delta }_h(0,0)}$ & $2$ 
& B \\
\rule{0pt}{11pt}$\gl (V, {\tau }_{\ast })_{v^0}$ & 
${\Delta }^{\varepsilon }_0(0)$ & $\underline{{\Delta }^{\varepsilon}_h(0), \, 
\lambda \in {\R }_{>0}} $ & $1$ & \\
\rule{0pt}{11pt} $\gl (V, {\tau }_{\ast })_{v^0}$ & ${\Delta }^{+}_0(0) + 
{\Delta }^{-}_0(0)$ & 
$\underline{{\Delta }^{+}_h(0) + {\Delta }^{-}_0(0)}$ & $2$ 
&  \\ %%%%%%%%%%%%%%%%%%%%%%
\rule{0pt}{11pt}$\oo (V, {\sigma }_{+}, \tau )_{v^0}$ & 
${\Delta }^{\varepsilon}_0(0)$ 
& $\underline{{\Delta }^{\varepsilon}_h(0), \,  \lambda \in {\R }_{>0}}, \, 
\mbox{$h$ even}$ & $1$ & A \\
\rule{0pt}{11pt}$\oo (V, {\sigma }_{+}, \tau )_{v^0}$ & ${\Delta }^{+}_0(0)+ 
{\Delta }^{-}_0(0)$ 
& $\underline{ {\Delta }^{+}_h(0)+ {\Delta }^{-}_h(0)}, \, \mbox{$h$ even} $& $2$ & 
A \\
\rule{0pt}{11pt}$\oo (V, {\sigma }_{-}, \tau )_{v^0}$ & 
${\Delta }_0(0,0)$ & $\underline{ {\Delta }_h(0,0)}, \,  \mbox{$h$ even} $ & 
$2$ & B \\
\rule{0pt}{11pt}$\oo (V, {\sigma }_{-}, \tau )_{v^0}$ & ${\Delta}^{\varepsilon }_0(0,0)$ 
& $\underline{ {\Delta }^{\varepsilon}_h(0,0)}, \mbox{$h$ odd} $ 
& $2$ & B\\ %%%%%%
\rule{0pt}{11pt}$\spp (V, {\sigma }_{+}, \tau )_{v^0}$ & 
${\Delta }^{\varepsilon}_0(0)$
& $\underline{{\Delta }^{\varepsilon}_h(0), \, \lambda \in {\R }_{>0}}, \, 
\mbox{$h$ odd}$ & $1$  & A \\
\rule{0pt}{11pt}$\spp (V, {\sigma }_{+}, \tau )_{v^0}$ & 
${\Delta }_0(0,0)$ & $\underline{{\Delta }_h(0,0)}, \, \mbox{$h$ even}$ & 
$2$ & A \\
\rule{0pt}{11pt}$\spp (V, {\sigma }_{-}, \tau )_{v^0}$ & 
${\Delta }_0(0,0)$ & $\underline{{\Delta }_h(0,0)}, \, \mbox{$h$ odd}$ & 
$2$ & A \\
\rule{0pt}{11pt}$\spp (V, {\sigma }_{-}, \tau )_{v^0}$ & 
${\Delta }^{\varepsilon}_0(0,0)$ & 
$\underline{{\Delta }^{\varepsilon }_h(0,0)}, \, \mbox{$h$ even}$ & 
$2$ & B 
\end{tabular}\label{pg-table2}\medskip 

\noindent \hspace{.2in}\parbox[t]{4.85in}{Table 2. List of indecomposable nilpotent distinguished types. Here $\varepsilon \in \R$ with ${\varepsilon }^2 =1$. Condition A is: ${\sigma }_{+}(v^0) = v^0$ and condition B is: $v^0 \in V_{{\sigma }_{-}}$ 
with ${\sigma }_{-}(v^0) = \pm \mathrm{i}v^0$ and $\tau (v^0, v^0) =0$.} \bigskip

\noindent \textsc{Case} $\mathbf{1}$. $\dim \overline{\overline{W}} =1$. \medskip 

\noindent Suppose that $\tau $ is not involved in the triple defining 
$\underline{\Delta }$ and $\sigma = {\sigma }_{+}$.   
\par Let $Y \in \gl (W; {\sigma }_{+})$ with $(0, \overline{\overline{W}}; 
{\overline{\overline{\sigma }}}_{+}) \in {\Delta}_0(0) $ and $(Y, W; {\sigma }_{+}) \in 
{\Delta }_h(0)$. The distinguished type $\underline{\Delta }$ has distinguished height $h$ and set of parameters $\mu (\underline{\Delta })$ with 
$v^0 \in  W_{{\sigma }_{+}} = \{ w \in W \setrule \, {\sigma }_{+}(w) = w \} $. 
The case where ${\overline{\sigma}}_{+}(v^0) = -v^0$ is handled by renaming 
the $-1$ eigenspace of ${\sigma }_{+}$ as $W_{{\sigma }_{+}}$. There is a nonzero vector $w \in W_{{\sigma }_{+}}$ and a nonzero $\lambda \in 
\mu (\underline{\Delta })$ such that $v_0 = \lambda Y^hw$. Since $\mathrm{dht}(\Delta ) = h$ and $Y \in \gl (W; {\sigma }_{+})$, it follows that 
$\widetilde{W} = {\spann }_{\C} \{ w, Yw, \ldots , Y^hw \} $ 
is a subspace of $W_{{\sigma }_{+}}$, which is $Y$ and ${\sigma }_{+}$ invariant. 
Moreover, $v^0 \in \widetilde{W}$. Since $\underline{\Delta }$ is 
indecomposable, $\widetilde{W} = W$. Now 
\begin{align*}
\overline{\lambda } Y^hw & = {\sigma }_{+}(\lambda Y^hw), \, \, 
\mbox{since ${\sigma }_{+}$ is anti-linear and $Y^h w \in W_{{\sigma }_{+}}$} \\
& = {\sigma }_{+}(v^0) = v^0, \, \, \mbox{because $(Y,W,v^0; {\sigma }_{+})$ is 
a triple} \\
& = \lambda Y^hw. 
\end{align*}
Thus $\overline{\lambda } = \lambda $, that is, $\lambda \in \R$. Denote 
$\underline{\Delta }$ by $\underline{{\Delta }_h(0), \, \lambda \ne 0}$. Here 
$\lambda $ is a real nonzero \emph{modulus}.  
\par The case $\gl (W; {\sigma }_{-})$ does not occur because the anti-linear 
mapping ${\overline{\overline{\sigma }}}_{-}$ is only defined when 
$\overline{\overline{W}}$ is even dimensional. \medskip 

Suppose that ${\tau }_{\ast}$ or $\tau $ occurs in the triple representing 
the distinguished type $\underline{\Delta }$ and that the induced bilinear form 
on $\overline{\overline{W}}$ is \emph{hermitian} or \emph{symmetric}, 
see appendix 1. The alternatives are \medskip %

\noindent \hspace{1in}\begin{tabular}{lcl}
\multicolumn{1}{c}{Lie algebra} & \multicolumn{1}{c}{
$\overline{\overline{\Delta}}$} & 
\multicolumn{1}{c}{\hspace{.2in}${\Delta }$} \\ \hline 
$\gl (W; {\tau }_{\ast } )$ & ${\Delta}^{\varepsilon}_0(0) $ & 
${\Delta }^{\varepsilon} _h(0)$   \\
$\oo (W; {\sigma }_{+}, \tau )$ & ${\Delta}^{\varepsilon}_0(0) $ & 
${\Delta }^{\varepsilon} _h(0)$, \, \mbox{$h$ even} \\
$\spp (W; {\sigma }_{+}, \tau )$ & ${\Delta}^{\varepsilon}_0(0) $ & 
${\Delta }^{\varepsilon} _h(0)$, \, \mbox{$h$ odd.} \\
\end{tabular}\medskip 

First we look at the case where $\underline{\Delta }$ is a distinguished type, represented by the triple $(Y,W,v^0; {\tau }_{\ast })$, of distinguished height 
$h$, where ${\tau }_{\ast }$ is hermitian. Hence $h$ is even. Then there is a nonzero vector $w \in W$ such that $Y^h w =v^0$. 
Thus the subspace $\widetilde{W} = {\spann }_{\C} \{ w, Yw, $ 
$\ldots , Y^hw = v^0 \}$ of $W$ is $Y$ invariant and contains the vector $v^0$. Since ${\overline{\overline{\tau}}}_{\ast }$ is nondegenerate on 
$\overline{\overline{W}}$ and $\dim \overline{\overline{W}} =1$, it follows that 
${\overline{\overline{\tau}}}_{\ast }(\overline{\overline{w}}, \overline{\overline{w}}) 
\in \R \setminus \{ 0 \}$, because ${\overline{\overline{\tau}}}_{\ast }$ is 
hermitian and $\overline{\overline{w}} \ne 0$ since $w \notin YW$. Thus we may assume that 
${\overline{\overline{\tau}}}_{\ast }(\overline{\overline{w}}, \overline{\overline{w}}) = 
\varepsilon \lambda $, where $\lambda > 0$ and ${\varepsilon }^2=1$. One has 
\begin{displaymath}
{\tau }_{\ast}(w,v^0) = {\tau }_{\ast}(w, Y^hw) = 
{\overline{\overline{\tau}}}_{\ast }(\overline{\overline{w}}, \overline{\overline{w}}) = 
\varepsilon \lambda , 
\end{displaymath}
which shows that $\varepsilon \lambda \in \mu (\underline{\Delta })$. Moreover, 
we have
\begin{displaymath}
{\tau }_{\ast }(Y^{i-1}w, Y^{j-1}w) = 
\mbox{\footnotesize $\left\{ \begin{array}{cl}
(-1)^{i-1}\varepsilon \lambda , & \mbox{if $(i-1)+(j-1) =h$} \\
0, & \mbox{otherwise}. \end{array}  \right. $}
\end{displaymath}
So the bilinear form ${\tau }_{\ast }$ on $\widetilde{W}$ is nondegenerate, 
since the matrix $T$ whose $ij^{\mathrm{th}}$ entry is 
$\big( ({\tau }_{\ast }(Y^{i-1}w, Y^{j-1}w) \big) $ has nonzero entries on its anti-diagonal and zero entries below its anti-diagonal, which implies $\det T \ne 0$. 
Hence $\widetilde{W} = W$, because 
$\underline{\Delta }$ is indecomposable. Denote $\underline{\Delta }$ by 
$\underline{{\Delta }^{\varepsilon}_h(0), \, \lambda >0 }$, $h$ is even. It has 
a positive modulus $\lambda $. Suppose that ${\overline{\overline{\tau }}}_{\ast }$ is alternating on $\overline{\overline{W}}$, 
which implies that $h$ is odd, since ${\tau}_{\ast}$ is 
hermitian. Then ${\widetilde{\tau }}_{\ast } = 
\mathrm{i}\, {\overline{\overline{\tau}}}_{\ast }$ is hermitian. Thus 
we may assume that ${\overline{\overline{\tau }}}_{\ast }$ is a nondegenerate hermitian form on $\overline{\overline{W}}$. Using the same argument as above, 
we find that the distinguished type $\underline{\Delta }$ is 
$\underline{{\Delta }^{\varepsilon}_h(0), \, \lambda >0}$, $h$ is odd. It has 
a positive modulus $\lambda $. \medskip 
 
The cases $\oo (W, v^0; {\sigma }_{+}, \tau )$ with $h$ odd and 
$\spp (W, v^0; {\sigma }_{+}, \tau )$ with $h$ even do not occur because 
$\overline{\overline{\tau }}$ is nondegenerate and alternating, 
which implies that $\dim \overline{\overline{W}}$ is even. 
Also the cases $\oo (W, v^0; {\sigma }_{-}, \tau )$ and $\spp (W, v^0; {\sigma }_{-}, 
\tau )$ do not occur, because $\dim \overline{\overline{W}}$ must be 
even for ${\overline{\overline{\sigma }}}_{-}$ to be defined. \medskip  

\noindent \textsc{Case} $\mathbf{2}$. 
$\dim \overline{\overline{W}} =2$.  \medskip 

\noindent Suppose $\tau $ is not involved in the triple defining the distinguished type $\underline{\Delta }$.  
\par We look at the case of $\gl (W; {\sigma }_{+})$. The reduced type 
$\overline{\overline{\underline{\Delta }}}$, represented by the pair 
$(\overline{\overline{Y}} =0, \overline{\overline{W}}; 
{\overline{\overline{\sigma }}}_{+})$ is not indecomposable. The following 
argument shows that $\underline{\Delta }$ is not indecomposable. Consequently, 
this case is excluded. The reduced type $\overline{\overline{\underline{\Delta }}}$ is the sum of two indecomposable types: ${\Delta }_0(0)$ and ${\Delta }_0(0)$. There 
are special vectors $\overline{\overline{z}}$ and $\overline{\overline{w}}$ in 
$\overline{\overline{W}}$ such that $\{ \overline{\overline{z}}, 
\overline{\overline{w}} \} $ is a basis of $\overline{\overline{W}}$. 
By \cite[prop.3 p.343]{burgoyne-cushman} there are $z, w \in W$ with 
$\overline{\overline{z}} = z +YW \in \overline{\overline{W}}$ and $\overline{\overline{w}} = w +YW \in \overline{\overline{W}}$ such that 
$\{ z, \, Yz, \ldots , Y^hz ; \, w, \, Yw, \ldots , Y^hw \} $ 
is a basis of $W$. Hence $Y^hz \ne 0$. Since $W = W\cap W^1 \oplus 
W \cap W^2$, where $W^i$ are eigenspaces of $\overline{\sigma }|W$, 
we may write $z = z_1+z_2$, where $z_i \in W\cap W^i$ for $i=1,2$. Now 
$0 \ne Y^hz = Y^hz_1 +Y^hz_2$ implies $Y^hz_i$ are not both zero. Suppose 
that $Y^hz_1 \ne 0$. Otherwise interchange $z_1$ and $z_2$ and rename 
$z_1$ as $z$. Then $z$ is a special vector. A similar argument shows that 
we may assume that $w$ is a special vector. We now show that we can choose a basis of $W$ so that $Y^hw = v^0$. Since $\ker Y = \spann \{ Y^hz, Y^hw \} $ and $v^0 \in \ker Y$ by hypothesis, there are $\alpha $, $\beta \in \C$ not both zero 
such that $v^0 = \alpha Y^hz + \beta Y^hw$. If $\alpha \ne 0$ let $z' = 
{\alpha }^{-1}z$ and $w' = \alpha w + \beta z$; while if $\alpha =0$ and 
$\beta \ne 0$ let $z' = \beta z$ and $w' = -{\beta }^{-1}z$. Then 
$W' = {\spann}_{\C}\{ Y^hz', \,  Y^{h-1}z', \ldots , z' ; \, w', \, Yw', \ldots , 
Y^hw' \} $ is a $Y$ and $\sigma $ invariant subspace of $W$, since 
$z'$ and $w'$ are special vectors. Moreover, $Y^hw'= v^0$. Since 
the distinguished type $\underline{\Delta }$ represented by the triple $(Y,W, v^0; {\sigma }_{+})$ is indecomposable, $\widetilde{W} = W$. But 
$\underline{\Delta }$ is the sum of the distinguished type $\underline{\Delta }'$, represented by the triple $(Y, W'', v^0; {\sigma }_{+})$, where $W'' = \spann \{ w', \, Yw', \ldots , Y^hw' = v^0 \}$, and the type $\widetilde{\Delta }$, represented by the 
pair $(Y, W'^{\dagger}; {\sigma }_{+})$, where $W^{\dagger} = 
{\spann}_{\C}\{ z', Yz', \ldots , Y^hz' \}$. In other words, $\underline{\Delta }$ is 
decomposable. 
\par We now look at the case $\gl (W; {\sigma }_{-}) $. The reduced type 
$\overline{\overline{\underline{\Delta }}}$, 
represented by the pair $(\overline{\overline{Y}} =0, \overline{\overline{W}}; 
{\overline{\overline{\sigma }}}_{-})$ is indecomposable. There is a basis 
$\{ \overline{\overline{z}}, \overline{\overline{w}} \} $ of special vectors 
of the vector space $\overline{\overline{W}}$ such that $\overline{\overline{z}} + { \overline{\overline{\sigma }} }_{-}(\overline{\overline{w}})$ is a basis of the vector space ${ \overline{\overline{W} }}_{ { \overline{\overline{\sigma }} }_{-} }$ 
over the quaternions, see appendix 1. By \cite[prop. 3 p.343]{burgoyne-cushman}
there are vectors $z$, $w \in W$ with $\overline{\overline{z}} = z +YW$ 
and $\overline{\overline{w}} = w +YW$ such that 
$\{ z, \, Yz, \ldots , Y^hz; \, w, \, Yw, \ldots , Y^hw \}$ is a basis of $W$. 
Because $Y^hz \ne 0$ and $Y^hw \ne 0$, we may 
assume that the vectors $z$ and $w$ are special. Arguing 
as in the preceding paragraph, we may choose $z$ and $w$ so that 
$Y^hw = v^0$. We denote the distinguished type $\underline{\Delta }$ by 
$\underline{{\Delta }_h(0,0)}$. There is no modulus. \medskip 
 
Suppose that $\tau $ is involved in the triple defining the distinguished type 
$\underline{\Delta }$ and that the bilinear form on $\overline{\overline{W}}$ is \emph{alternating}.  
\par Suppose that ${\tau }_{\ast }$ is hermitian. Then $h$ is odd, because 
$\overline{\overline{\tau}}$ is alternating. Hence ${\widetilde{\tau }}_{\ast } = 
\mathrm{i}\, \overline{\overline{\tau }}_{\ast }$ is a nondegenerate symmetric hermitian form on $\overline{\overline{W}}$, where the reduced type 
$\overline{\overline{\Delta }}$ is represented by the triple $(\overline{\overline{W}}; 
\overline{\overline{\sigma }}, \overline{\overline{\tau }})$.
Using the same argument as in case 1, we find that distinguished type 
$\underline{\Delta }$ is the sum of an indecomposable distinguished type 
${\Delta }^{{\varepsilon}_1}_h(0), \, {\lambda}_1 >0$ and an indecomposable 
type ${\Delta }^{{\varepsilon }_2}_h(0)$. Hence $\underline{\Delta }$ is not indecomposable. This contradicts our hypothesis. Thus 
the case $\gl (W, {\tau }_{\ast })$, $h$ odd does not occur. \medskip

The remaining alternatives are listed below. \medskip %

\noindent \hspace{1in}\begin{tabular}{lcl}
\multicolumn{1}{c}{Lie algebra} & \multicolumn{1}{c}{
$\overline{\overline{\Delta}}$} & 
\multicolumn{1}{l}{$\hspace{.5in}{\Delta }$} \\ \hline 
$\oo (W; {\sigma }_{+}, \tau )$ & ${\Delta }_{0}(0,0)$ & ${\Delta }_h(0,0)$, \, 
\mbox{$h$ odd} \\
$\oo (W; {\sigma }_{-}, \tau )$ & ${\Delta}^{\varepsilon}_0(0, 0) $ & 
${\Delta }^{\varepsilon}_h(0, 0)$, \, \mbox{$h$ odd} \\
$\spp ( W; {\sigma }_{+}, \tau )$ & ${\Delta}_0(0,0) $ & 
${\Delta }_h(0,0)$, \, \mbox{$h$ even} \\
$\spp ( W; {\sigma }_{-}, \tau )$ & ${\Delta}^{\varepsilon}_0(0,0) $ & 
${\Delta }^{\varepsilon}_h(0,0)$, \, \mbox{$h$ even.} \\
\end{tabular}\medskip

\noindent Suppose that $\sigma = {\sigma }_{+}$ and that the reduced type 
$\overline{\overline{\underline{\Delta }}}$ 
is represented by the pair $(\overline{\overline{Y}} =0, \overline{\overline{W}}; 
{\overline{\overline{\sigma }}}_{+},$ $\overline{\overline{\tau }})$. 
There are nonzero vectors $\overline{\overline{z}}$ and $\overline{\overline{w}}$ 
in $\overline{\overline{W}}$ such that $\{ \overline{\overline{z}}, 
\overline{\overline{w}} \} $ is a  
basis for ${\overline{\overline{W}}}_{{\sigma }_{+}}$ with $\overline{\overline{\tau }}(\overline{\overline{z}}, \overline{\overline{z}}) = 0 = 
\overline{\overline{\tau }}(\overline{\overline{w}}, \overline{\overline{w}})$ and 
$\overline{\overline{\tau }}(\overline{\overline{z}}, \overline{\overline{w}}) = \mu \ne 0$. By \cite[prop.3 p.343]{burgoyne-cushman} there are vectors $z$, $w \in 
W$ so that $\overline{\overline{z}} = z +YW \in \overline{\overline{W}}$ and 
$\overline{\overline{w}} = w +YW \in \overline{\overline{W}}$ such that 
$\widetilde{W} = {\spann }_{\C} \{ z, \, Yz, \ldots , Y^hz; \, w, \, Yw, \ldots , Y^hw \} $ 
is a subspace of $W$. Since $Y^hz$ and $Y^hw$ are both nonzero, we 
may assume that $z$, $w \in W_{{\sigma }_{+}}$. Moreover,  
\begin{displaymath}
\tau (Y^{i-1}z, Y^{j-1}w) = \mbox{\footnotesize $\left\{ \begin{array}{rl}
\hspace{-5pt} (-1)^i  \mu , &\mbox{if $(i-1)+(j-1) =h$} \\ 0, & \mbox{if $i+j > h$,} \end{array} \right. $} 
\end{displaymath}
$\tau (Y^{i-1}z,Y^{j-1}z) = 0 = \tau (Y^{i-1}w, Y^{j-1}w)$. We now show that we can choose a basis of $W_{{\sigma }_{+}}$ so that $Y^hw = v^0$. There are $\alpha $, $\beta \in \C$ not both zero such that 
$v^0 = \beta Y^hz + \alpha Y^hw$, since 
$v^0 \in \ker Y^h$ and $\ker Y^h = \spann \{ Y^hz, Y^hw \} $. Let  
\raisebox{1pt}{{\tiny $\begin{pmatrix} z' \\ w'' \end{pmatrix}$}$ =
${\tiny$\begin{pmatrix} a & b \\ c & d \end{pmatrix} $}\hspace{-4pt} 
{\tiny $\begin{pmatrix} z \\ w \end{pmatrix}$}}, where 
\raisebox{1pt}{{\tiny $\begin{pmatrix} a & b \\ c & d \end{pmatrix}$}$=$ 
{\tiny $\begin{pmatrix} {\alpha }^{-1} & 0 \\ \beta & \alpha \end{pmatrix}$}} if 
$\alpha \ne 0$ and {\tiny $\begin{pmatrix} \beta & 0 \\ 0 & {\beta }^{-1} 
\end{pmatrix}$} if $\beta \ne 0$. Then $\{ z', \, Yz', \ldots , Y^hz'; \, w',$ $ \, Yw', \ldots ,$ $ Y^hw' \} $ 
is spans $\widetilde{W}$ because $\tau (Y^{i-1}z', Y^{j-1}z') $ $= 0 = 
\tau (Y^{i-1}w', Y^{j-1}w') $ and 
\begin{align*}
\tau (Y^jz', Y^{h-j}w') & = (-1)^j \tau (a z + bw, Y^{h}(c z + d w)) \\ 
& = ac \, \overline{\overline{\tau }}(\overline{\overline{z}}, \overline{\overline{z}}) 
+bd \,  \overline{\overline{\tau }}(\overline{\overline{w}}, \overline{\overline{w}}) 
+ (ad-bc) \overline{\overline{\tau }}(\overline{\overline{z}}, \overline{\overline{w}}) \\
& = \overline{\overline{\tau }}(\overline{\overline{z}}, \overline{\overline{w}}), \, \, 
\parbox[t]{3in}{since $\overline{\overline{\tau }}$ is alternating and 
$ad -bc =1$} \\
& = \mu ;  
\end{align*}
while $\tau (Y^{i-1}z', y^{j-1}w') = 0$ if $(i-1)+(j-1) > h$. Thus the matrix 
$T =${\tiny $\left( \begin{array}{c|c}
\tau (Y^{i_1}z', Y^{j-1}z') & \tau (Y^{i-1}z', Y^{j-1}w') \\ \hline 
\tau (Y^{i-1}w', Y^{j-1}z') & \tau (Y^{i-1}w', Y^{j-1}w') \end{array} \right)$} on 
$\widetilde{W}$ is 
\begin{displaymath}
\mbox{\footnotesize $\left( \begin{array}{ccc|ccc}
& & & \mbox{\Large $\ast $} & & + \\
& \mbox{\Large $0$} & &  &+ &  \\
& & & + & & \mbox{\Large $0$} \\ \hline 
\mbox{\Large $\ast $} & & + & & &  \\
& + & & & \mbox{\Large $0$} & \\
+ & & \mbox{\Large $0$} & & & 
\end{array} \right) $. }
\end{displaymath}
So $\det T \ne 0$, since its antidiagonal entries are all nonzero. Thus $\tau $ is nondegenerate on $\widetilde{W}$. By construction $v^0 = Y^hw'$. But 
$\widetilde{W}$ is a $Y$ invariant subspace of $W$, which contains the vector $v^0$. It is also ${\sigma }_{+}$ invariant, since the vectors $z'$ and $w'$ lie in 
$W_{{\sigma }_{+}}$ and thus are special. Hence $\widetilde{W} = W$, since 
$\underline{\Delta }$ is indecomposable. We denote $\underline{\Delta }$ by 
$\underline{{\Delta }_h(0,0)}$, when $h$ is odd for $\oo (W, v^0; {\sigma }_{+}, 
\tau )$ or when $h$ is even for $\spp (W, v^0; {\sigma }_{+}, \tau )$. There is no modulus. \medskip 

Suppose that $\sigma = {\sigma }_{-}$ and that the reduced type 
$\overline{\overline{\underline{\Delta }}}$ is 
represented by the pair $(\overline{\overline{Y}}=0, \overline{\overline{W}}; $ $ 
{\overline{\overline{\sigma }}}_{-}, \overline{\overline{\tau }})$, where 
$\overline{\overline{\tau}}$ is alternating. There is a basis 
$\{ \overline{\overline{z}}, \overline{\overline{w}} \} $ of eigenvectors of 
${\sigma }_{-}$ on $\overline{\overline{W}}$ such that $\overline{\overline{y}} = 
\overline{\overline{z}} +{\overline{\overline{\sigma }}}_{-}(\overline{\overline{w}})j$ 
is a basis for the $1$ dimensional quaternionic vector space 
${\overline{\overline{W}}}_{{\sigma }_{-}}$ 
with a nondegenerate hamiltonian alternating form $\boldsymbol{\tau }$, see appendix $1$. Thus 
$\boldsymbol{\tau}(\overline{\overline{y}}, \overline{\overline{y}}) = 
\lambda \varepsilon  j$, where $\lambda \in {\R }_{>0}$ and 
${\varepsilon}^2 =1$. Let ${\overline{\overline{z}}}\, ' = 
{\lambda }^{-1/2}\, \overline{\overline{z}}$ and 
${\overline{\overline{w}}}\, ' = {\lambda }^{-1/2}\, \overline{\overline{w}}$, where
$\overline{\overline{z}}$ and $\overline{\overline{w}}$ are special vectors. Then with $y' = \overline{\overline{z}}\, ' + 
{\overline{\overline{\sigma }}}_{-}(\overline{\overline{w}}\, ') j$ one has 
$\boldsymbol{\tau}(\overline{\overline{y}}\, ', \overline{\overline{y}}\, ') = 
\varepsilon  j$. Since $\tau (z', Y^hw') = 
\overline{\overline{\tau}}({\overline{\overline{z}}}\, ', {\overline{\overline{w}}}\, ') 
\ne 0$ for $z' \in {\overline{\overline{z}}}\, '$ and $w' \in {\overline{\overline{w}}}\, '$ 
with $Y^hz'$ and $Y^hw$ both nonzero, an argument shows that we may assume that $z'$ and $w'$ are special vectors in $W$ and that 
$W' = {\spann }_{\C}\{ z', \ldots , Y^hz'; w, \ldots , Y^hw' \} $ contains the 
vector $v^0$. Since the distinguished type $\underline{\Delta }$, 
represented by the triple $(Y, W, v^0; \sigma , \tau )$ is indecomposable, 
it follows that $W' = W$ and there is no modulus. Hence the distinguished type 
$\underline{\Delta }$ is $\underline{{\Delta }^{\varepsilon}_h(0,0)}$ with 
$h$ odd, when $\tau $ is symmetric and $h$ is even when $\tau $ is 
altenating. \medskip 

Suppose that $\tau $ or ${\tau }_{\ast }$ is involved in the triple defining the distinguished type $\underline{\Delta }$ and that the 
bilinear form on $\overline{\overline{W}}$ is \emph{symmetric} or 
\emph{hermitian}. The alternatives are listed below.  \medskip 

\noindent \hspace{.5in}\begin{tabular}{lll}
\multicolumn{1}{c}{Lie algebra} & \multicolumn{1}{c}{
$\overline{\overline{\Delta}}$} & 
\multicolumn{1}{c}{${\Delta }$} \\ \hline 
$\gl (W; {\tau }_{\ast })$ & ${\Delta }^{+}_0(0)+{\Delta }^{-}_0(0)$ & 
${\Delta }^{+}_h(0)+{\Delta }^{-}_h(0)$, \, \mbox{$h$ even} \\
$\oo (W; {\sigma }_{+}, \tau )$ & ${\Delta}^{+}_0(0) +{\Delta }^{-}_0(0) $ & 
${\Delta }^{+}_h(0) +{\Delta }^{-}_h(0)$, \, \mbox{$h$ even} \\
$\oo (W; {\sigma }_{-}, \tau )$ & ${\Delta}_0(0,0) $ & 
${\Delta }_h(0,0)$, \, \mbox{$h$ even} \\
$\spp (W; {\sigma }_{+}, \tau )$ & ${\Delta}^{+}_0(0) +{\Delta }^{-}_0(0) $ & 
${\Delta }^{+}_h(0) +{\Delta }^{-}_h(0)$, \, \mbox{$h$ odd}  \\
\rule{0pt}{10pt}$\spp (W; {\sigma }_{-}, \tau )$ & ${\Delta }_0(0,0)$ & 
${\Delta }_h(0,0)$, \, \mbox{$h$ odd.}
\end{tabular}\medskip

\noindent We treat the case when the bilinear form on 
$\overline{\overline{W}}$ is the hermitian form ${\tau }_{\ast }$. 
Suppose that ${\overline{\overline{\tau }}}_{\ast}(\overline{\overline{z}}, 
\overline{\overline{z}}) \ne 0$ for every nonzero $\overline{\overline{z}} \in \overline{\overline{W}}$. Then there is a basis $\{ \overline{\overline{z}}, 
\overline{\overline{w}} \} $ of $\overline{\overline{W}}$, which is 
${\overline{\overline{\tau }}}_{\ast}$ orthogonal, that is, 
${\overline{\overline{\tau }}}_{\ast}(\overline{\overline{z}}, 
\overline{\overline{z}}) = \mu \ne 0$, ${\overline{\overline{\tau }}}_{\ast}(\overline{\overline{w}}, \overline{\overline{w}}) = \lambda \ne 0$, and 
${\overline{\overline{\tau }}}_{\ast}(\overline{\overline{z}}, 
\overline{\overline{w}}) =0$. Then there 
are vectors $z$, $w \in W$ with $\overline{\overline{z}} = z +YW$ and 
$\overline{\overline{w}} = w +YW$ such that $\{ z, Yz, \ldots , Y^hz; 
w, Yw, \ldots , Y^hw \}$ forms a ${\tau }_{\ast }$ orthogonal basis 
of $W$. We can choose $\overline{\overline{w}}$ so that $Y^hw = v^0$. 
Thus $\underline{\Delta }$ is the sum of the distinguished type 
${\underline{\Delta }}_1$, represented by the triple 
$(Y|W_1, W_1, v^0; {\tau }_{\ast }|W_1)$, where $W_1 = 
\spann \{ w, Yw, \ldots , Y^hw = v^0 \}$, and a type ${\Delta }_2$, 
represented by the pair $(Y|W_2, W_2; {\tau }_{\ast }|W_2)$, where 
$W_2 = \spann \{ z, Yz, \ldots , Y^hz \}$. This contradicts our hypothesis 
that $\underline{\Delta }$ is indecomposable. Thus there is a nonzero 
vector $\overline{\overline{z}} \in \overline{\overline{W}}$ such that 
${\overline{\overline{\tau }}}_{\ast}(\overline{\overline{z}}, 
\overline{\overline{z}}) =0$. Since ${\overline{\overline{\tau }}}_{\ast}$ is 
nondegenerate on $\overline{\overline{W}}$, there is a vector 
$\overline{\overline{y}} \in \overline{\overline{W}}$ such that 
${\overline{\overline{\tau }}}_{\ast}(\overline{\overline{z}}, 
\overline{\overline{y}}) = \eta \ne 0$. Let $\overline{\overline{w}} = 
\overline{\overline{y}}  - \frac{ {\overline{\overline{\tau }}}_{\ast } (\overline{\overline{y}}, \overline{\overline{y}}) }{\rule{0pt}{8pt} 2 \, {\overline{\overline{\tau }}}_{\ast }(\overline{\overline{z}}, \overline{\overline{y}})} \, \overline{\overline{z}}$. 
Then $\overline{\overline{w}}$ is a ${\overline{\overline{\tau }}}_{\ast }$ isotropic vector in $\overline{\overline{W}}$. Thus the matrix of 
${\overline{\overline{\tau }}}_{\ast }$ with respect to the basis 
$\{ \overline{\overline{z}}, \overline{\overline{w}} \} $ is \raisebox{1pt}{{\tiny $\begin{pmatrix} 
0 & \eta \\ \eta & 0 \end{pmatrix}$}}, because ${\overline{\overline{\tau }}}_{\ast }
(\overline{\overline{z}}, \overline{\overline{w}}) = {\overline{\overline{\tau }}}_{\ast }
(\overline{\overline{z}}, \overline{\overline{y}} ) = \eta $. Hence the semisimple reduced type $\overline{\overline{\underline{\Delta }}}$ is equal to 
${\Delta }^{+}_0(0) + {\Delta }^{-}_0(0)$.  Using \cite[prop.2 p.343]{burgoyne-cushman} there are vectors $z$, $w \in W$ with $\overline{\overline{z}} = z + YW \in 
\overline{\overline{W}}$ and $\overline{\overline{w}} = w + YW \in 
\overline{\overline{W}}$ such that 
$\{ Y^hz, \, Y^{h-1}z, \ldots , z; \, w, \, Yw, \ldots , Y^hw \}$ 
is a basis of $W$ where ${\tau }_{\ast }(Y^iz, Y^jw)=
$\raisebox{1pt}{{\tiny $\left\{ 
\begin{array}{cl} \hspace{-5pt}(-1)^i \eta  , & \hspace{-5pt} \mbox{if $i+j =h$} \\ 
\hspace{-5pt} 0, & \hspace{-5pt} \mbox{otherwise} \end{array} 
\right. $}} and ${\tau }_{\ast }(Y^iz, Y^j z ) = 0 = {\tau }_{\ast } (Y^iw, Y^jw)$. Since 
$v^0 \in \ker Y = \spann \{ Y^hz, Y^hw \} $, there are $\alpha $, $\beta \in \C$ 
not both zero such that $v^0 = \alpha Y^hz + \beta Y^hw$. We now show that 
we can find a basis of $W$ such that $v^0 = Y^hw$. If $\beta 
\ne 0$ set $z' = {\beta }^{-1}z$ and $w' =\alpha z + \beta w$; while 
if $\alpha \ne 0$ and $\beta =0$, set $w' = \alpha z$ and $z' = -{\alpha }^{-1}w$. 
In either case $v^0 = Y^hw'$ and 
$\{ Y^hz', \, -Y^{h-1}z', \ldots , z'; \, w', \, Yw', \ldots , Y^hw' \} $ 
is a basis of $W$, which is $Y$ invariant, and with respect to which the matrix of 
${\tau }_{\ast }$ is 
\raisebox{1pt}{{\tiny $\begin{pmatrix} 0 & \eta I_{h+1} \\ \eta  I_{h+1} & 0 \end{pmatrix}$}}. 
The indecomposable nilpotent distinguished type $\underline{\Delta }$ is denoted 
$\underline{ {\Delta }^{+}_h(0) + {\Delta }^{-}_h(0) }$. There is no modulus. \medskip 

Suppose that $\overline{\overline{\tau }}$ is symmetric and 
$\sigma = {\sigma }_{+}$. In the preceding paragraph replace 
${\overline{\overline{\tau }}}_{\ast }$ by ${\overline{\overline{\tau }}}_{+} = 
\overline{\overline{\tau }}|W_{{\sigma }_{+}}$, ${\tau }_{\ast }$ by 
${\tau }_{+}$, and $W$ by $W_{{\sigma }_{+}}$. Arguing as before, we find 
that $\overline{\overline{\underline{\Delta }}}$ is 
$\underline{ {\Delta }^{+}_h(0) + {\Delta }^{-}_h(0) }$, when $h$ is even 
and $\tau $ is symmetric or $\tau $ is alternating and $h$ is odd. There is no modulus. \medskip 

Suppose that $\overline{\overline{\tau }}$ is a symmetric bilinear form 
on $\overline{\overline{W}}$ and $\sigma = {\sigma }_{-}$. 
Here the triple $(W, v^0; {\sigma }_{-}, \tau )$ represents the 
indecomposable nilpotent distinguished type $\underline{\Delta }$ of distinguished height $h$, where $h$ even when $\tau $ is symmetric or $h$ is odd when 
$\tau $ is alternating. The pair $(\overline{\overline{Y}} =0, \overline{\overline{W}}; {\overline{\overline{\sigma }}}_{-}, \overline{\overline{\tau }} )$ represents the reduced type $\overline{\overline{\underline{\Delta }}}$. There is a basis 
$\{ \overline{\overline{z}}, \overline{\overline{w}} \} $ of special vectors of 
$\overline{\overline{W}}$ such that $\overline{\overline{y}} = 
\overline{\overline{z}} + {\overline{\overline{\sigma }}}_{-}(\overline{\overline{w}})j$ 
is a basis of the $1$ dimensional quaternionic vector space 
${\overline{\overline{W}}}_{{\sigma }_{-}}$ with hamiltonian symmetric form 
$\boldsymbol{\tau}$ such that $\boldsymbol{\tau}(\overline{\overline{y}}, 
\overline{\overline{y}}) = \lambda $, where $\lambda \in {\R }_{>0}$. 
Let $\overline{\overline{z}}\, ' = {\lambda }^{-1/2}\, \overline{\overline{z}}$, 
$\overline{\overline{w}}\, ' = {\lambda }^{-1/2}\, \overline{\overline{w}}$, 
and $\overline{\overline{y}}\, ' = {\lambda }^{-1/2}\overline{\overline{y}}$. 
Then $\overline{\overline{\tau }}(\overline{\overline{z}}\, ', 
\overline{\overline{z}}\, ') = 1 = \overline{\overline{\tau }}(\overline{\overline{w}}\, ', 
\overline{\overline{w}})$ and $\overline{\overline{\tau }}(\overline{\overline{z}}, 
\overline{\overline{w}}\, ' ) =0$ and the reduced type 
$\overline{\overline{\underline{\Delta }}}$ is ${\Delta }_0(0,0)$. 
Using \cite[prop.2 p.343]{burgoyne-cushman}, we find that 
there are vectors $z$, $w \in W$ with $\overline{\overline{z}} = z +YW$ and 
$\overline{\overline{w}} = w +YW$ such that $\{ Y^hz, \, Y^{h-1}z, \ldots , z; 
w, \, Yw, \ldots , Y^hw \} $ is a basis of $W$ such that 
$\tau (Y^iz,Y^jw)=${\tiny $\left\{ \begin{array}{cl} 
(-1)^i \lambda , & \mbox{if $i+j = h$} \\ 0, & \mbox{otherwise} \end{array} \right. $}, 
where $\lambda \in {\R }\setminus \{ 0 \}$, and 
$\tau (Y^iz, Y^jz) = 0 = \tau (Y^iw, Y^jw)$. Since $\tau (z, Y^hw) \ne 0$, we may assume that $z$ and $w$ are special vectors. Then 
$W$ is ${\sigma }_{-}$ invariant. We now show that we can choose 
a $\tau $ orthogonal basis of $W$ such that $Y^hw = v^0$. 
Since $v^0 \in \ker Y^h = \spann \{ Y^hz , \, Y^hw \} $, there are $\alpha $, 
$\beta \in \C$ not both zero such that $v^0 = \alpha Y^hz + \beta Y^hw$. 
If $\beta \ne 0$, set $z'= {\beta }^{-1}z$ and $w' =\alpha z + \beta w$. 
If $\alpha \ne 0$ and $\beta =0$, set $z' = -{\alpha }^{-1}w$ and $w' = 
\alpha z$. In either case $Y^hw' = v^0$ and $\{ Y^hz', -Y^{h-1}z', \ldots , z'; 
w', \, Yw' , \ldots , Y^hw' \} $ is a basis of $W$ such that the matrix 
of $\tau $ is {\tiny $\begin{pmatrix} 0 & \lambda I_{m+1} \\ 
\lambda I_{m+1} & 0 \end{pmatrix}$}. 
The nilpotent distinguished type $\underline{\Delta }$ is denoted 
$\underline{{\Delta }_h(0,0)}$. There is no modulus.  \medskip 

This completes the classification of nilpotent indecomposable distinguish-ed 
special types and proves \medskip 

\noindent \textbf{Proposition 3.2} An indecomposable nilpotent distinguished 
type $\underline{\Delta }$ is uni-quely determined by its distinguished height 
$h$, an element of its set of parameters $\mu (\underline{\Delta })$, if 
nonempty, and the dimension of its reduced type 
$\overline{\overline{\underline{\Delta }}}$. \medskip 

\noindent \textbf{Proof of theorem 2.1} We prove the uniqueness of the 
decomposition of a distinguished type $\underline{\Delta }$ into a 
sum of an indecomposable nilpotent distinguished type ${\underline{\Delta }}\, '$ 
and a type $\Delta $, which is the sum of indecomposable types 
$\sum^k_{\ell = 1}{\Delta }_{\ell }$. From the 
fact that ${\underline{\Delta }}'$ has the same distinguished type and 
set of parameters as $\underline{\Delta }$, using proposition 3.2 it follows 
that ${\underline{\Delta }}'$ is unique. From the theorem 
\cite[p.343]{burgoyne-cushman} it follows that the decomposition of 
the type $\Delta $ into a sum of indecomposable types is unique up to 
reordering of the summands. Thus the decomposition 
$\underline{\Delta } = {\underline{\Delta }}' + \sum^k_{\ell =1}{\Delta }_{\ell }$ is 
unique. \hfill $\square $

\section{The real affine classical groups}
%%%%%%%%%%%%%%%%%%

Let $\widetilde{V}$ be a complex vector space of dimension $n$ with an 
anti-linear mapping $\widetilde{\sigma }$ of $\widetilde{V}$ into itself 
such that ${\widetilde{\sigma }}^{\, 2} = \pm {\mathrm{id}}_{\widetilde{V}}$ and 
a nondegenerate bilinear form $\tau : \widetilde{V} \times \widetilde{V} 
\rightarrow \C $ that is symmetric or alternating and satisfies 
${\widetilde{\sigma }}^{\, \ast}\tau = \overline{(\tau )}$. $\widetilde{V}$ 
may also have a nondegenerate hermitian form ${\widetilde{\tau }}_{\ast }$, 
see appendix 1. \medskip 

Let $G(\widetilde{V}, \widetilde{\sigma }, \widetilde{\tau })$ be a 
real classical group, see table 1. In other words, $g \in 
G(\widetilde{V}, \widetilde{\sigma }, \widetilde{\tau })$ if and only if 
$g \in \Gl (\widetilde{V})$ such that ${\widetilde{\sigma }}^{\, \ast}(g) = g$ and 
$g^{\ast  } \widetilde{\tau } = \widetilde{\tau }$. The \emph{affine 
real classical group} $\mathrm{Aff}G(\widetilde{V}, \widetilde{\sigma }, \widetilde{\tau })$ is the set of complex affine mappings 
$\mathrm{aff}g: \widetilde{V} \rightarrow \widetilde{V}: 
\widetilde{v} \mapsto g(\widetilde{v}) + \widetilde{u}$, 
where $g \in G(\widetilde{V}, \widetilde{\sigma }, \widetilde{\tau })$ and 
$\widetilde{u}$, $\widetilde{v} \in \widetilde{V}$, with 
group multiplication given by composition of affine mappings. For $v^0$ a 
special vector in $\widetilde{V}$ the isotropy group ${G(\widetilde{V}, 
\widetilde{\sigma }, \widetilde{\tau})}_{v^0}$ is the set of all 
$g \in G(\widetilde{V}, \widetilde{\sigma }, \widetilde{\tau})$ such that 
$g(v^0) = v^0$. \medskip 

In this section we prove \medskip 

\noindent \textbf{Theorem 4.1}  Every real affine classical group 
$\mathrm{Aff}G(\widetilde{V}, \widetilde{\sigma }, \widetilde{\tau })$ is isomorphic 
to an isotropy group of a real classical group. \medskip

\noindent The proof procedes case by case. \medskip  

\noindent \textsc{Case} $\mathbf{1}$. $\mathrm{Aff}G(\widetilde{V}, {\widetilde{\sigma }}_{+})$. 
Let $V^{\vvee} = \widetilde{V} \times \C $. Set 
\begin{displaymath}
{\sigma}^{\vvee}_{+}: V^{\vvee} \rightarrow V^{\vvee}: v^{\vvee} = 
(\widetilde{v}, z) \mapsto {\sigma }^{\vvee}_{+}(v^{\vvee}) = 
({\widetilde{\sigma }}_{+}(\widetilde{v}), \overline{z}), 
\end{displaymath}
where 
\begin{displaymath}
{\widetilde{\sigma }}_{+}: \widetilde{V} \rightarrow \widetilde{V}: 
(z_1, \ldots , z_n)^T \mapsto ({\overline{z}}_1, \ldots , {\overline{z}}_n)^T, 
\end{displaymath}
see appendix 1. Then ${\sigma}^{\vvee}_{+}$ is an anti-linear mapping of 
$V^{\vvee}$ into itself with 
$({\sigma }^{\vvee}_{+})^2 = {\mathrm{id}}_{V^{\vvee}}$. 
Let $V^{\vvee}_{{\sigma}^{\vvee}_{+}} = \{ v^{\vvee} \in V^{\vvee} \setrule \, 
{\sigma}^{\vvee}_{+}(v^{\vvee}) = v^{\vvee} \} $. Then 
$V^{\vvee}_{{\sigma}^{\vvee}_{+}}$ is an $n+1$ dimensional real vector space 
with basis ${\epsilon }^{\vvee} = \{ e_1, \ldots , e_n; \, e_{n+1} \}$, where 
$\widetilde{\epsilon } = \{ e_1, \ldots , e_n \} $ is a basis of $\widetilde{V}$. 
For a nonzero vector $(v^{\vvee})^0 \in V^{\vvee}_{{\sigma }^{\vvee}_{+}}$, 
which is special by definition, let 
${G(V^{\vvee}, {\sigma }^{\vvee}_{+})}_{(v^{\vvee})^0}$ is 
the set of all $g \in \Gl (V^{\vvee})$ such that ${\sigma }^{\vvee}_{+}g = 
g{\sigma }^{\vvee}_{+}$ and $g((v^{\vvee})^0) = (v^{\vvee})^0$. Choose 
$(v^{\vvee})^0 = e_{n+1}$. Then $(v^{\vvee})^0$ is a special vector and 
\begin{displaymath}
{G(V^{\vvee}, {\sigma}^{\vvee}_{+})}_{(v^{\vvee})^0}  = 
\Big\{ 
\mbox{\footnotesize $\begin{pmatrix} \widetilde{A} & 0 \\ {\widetilde{d}}^{\, T} & 1 \end{pmatrix}$} 
\in \Gl (V^{\vvee}_{{\sigma}^{\vvee}_{+}}) \setrule \, 
\widetilde{A} \in \Gl ({\widetilde{V}}_{{\widetilde{\sigma }}_{+}} ) \, \, \& \, \, 
\widetilde{d} \in {\widetilde{V}}_{{\sigma}^{\vvee}_{+}} \Big\} 
\end{displaymath}
is isomorphic to $\mathrm{Aff}G(\widetilde{V}, {\widetilde{\sigma }}_{+})$. \medskip 

\noindent \textsc{Case} $\mathbf{2}$. $\mathrm{Aff}G(\widetilde{V}, {\widetilde{\sigma }}_{+}, 
\widetilde{\tau })$. Let $\widehat{V} = \C \times \widetilde{V} \times \C$. 
For $0 \le p \le \left[ \frac{n}{2} \right]$ set  
\begin{displaymath}
{\widehat{\sigma}}_{+} = {\widehat{\sigma}}^{(p)}_{+}: 
\widehat{V} \rightarrow \widehat{V}: \widehat{v} = 
(z_0, \widetilde{v}, z_{n+1})^T \mapsto {\widehat{\sigma}}_{+}(\widehat{v}) = 
({\overline{z}}_0, {\widetilde{\sigma }}^{(p)}_{+}(\widetilde{v}), {\overline{z}}_{n+1})^T,
\end{displaymath}
where 
\begin{displaymath}
{\widetilde{\sigma }}_{+} = {\widetilde{\sigma }}^{(p)}_{+}: \widetilde{V} \rightarrow \widetilde{V}: 
\widetilde{v} = (z_1, \ldots , z_n)^T \rightarrow  \mbox{ 
\footnotesize $\begin{pmatrix} -I_{n-p} & 0 \\ 0 & I_p \end{pmatrix}  
\begin{pmatrix} {\overline{z}}_1 \\ \vdots \\ {\overline{z}}_n \end{pmatrix} $,}
\end{displaymath}
see appendix 1. Then ${\widehat{\sigma}}_{+}$ is an anti-linear mapping of 
$\widehat{V}$ into itself such that $({\widehat{\sigma }}_{+})^2 = 
{\mathrm{id}}_{\widehat{V}}$. Set ${\widehat{V}}_{{\widehat{\sigma}}_{+}} = 
\{ \widehat{v} \in \widehat{V} \setrule \, {\widehat{\sigma}}_{+}(\widehat{v}) = \widehat{v} \} $. Then ${\widehat{V}}_{{\widehat{\sigma}}_{+}}$ is an $n+2$ dimensional real vector space. Let $\tau : \widehat{V} \times \widehat{V} \rightarrow \C$ be the complex valued bilinear form given by 
\begin{displaymath}
\tau (\widehat{u}, \widehat{v}) = \mbox{\footnotesize $\left\{ \begin{array}{rl}
w_0z_0 +\widetilde{\tau }(\widetilde{u}, \widetilde{v}) + w_{n+1}z_{n+1}, & 
\mbox{\normalsize if $\widetilde{\tau }$ is symmetric} \\
\rule{0pt}{10pt}w_0z_{n+1} +\widetilde{\tau }(\widetilde{u}, \widetilde{v})
-w_{n+1}z_0, & 
\mbox{\normalsize if $\widetilde{\tau }$ is alternating.}  \end{array} \right. $}
\end{displaymath}
The bilinear form $\tau $ is nondegenerate and symmetric if 
$\widetilde{\tau }$ is symmetric and alternating if $\widetilde{\tau }$ is alternating. It is easy to check that $\widehat{\tau }({\widehat{\sigma }}_{+}(\widehat{w}), 
{\widehat{\sigma }}_{+}(\widehat{v}))$ $= \overline{\widehat{\tau }(\widehat{w}), \widehat{v}))}$. \medskip 

Suppose that $\tau $ is \emph{symmetric}. Let ${\tau }_{+} = 
\tau |{\widehat{V}}_{{\widehat{\sigma }}_{+}}$, where 
${\widehat{V}}_{{\widehat{\sigma }}_{+}} = \{ \widehat{v} \in \widehat{V} \setrule \, 
{\widehat{\sigma }}_{+}(\widehat{v}) = \widehat{v} \}$. Then ${\tau }_{+}$ is a 
nondegenerate real valued symmetric bilinear form on 
${\widehat{V}}_{{\widehat{\sigma }}_{+}}$. 
For $0 \le p \le \left[ \frac{n}{2} \right] $ suppose that ${\widetilde{\tau }}_{+} = {\widetilde{\tau }}^{\, (p)}_{+}$ is a nondegenerate, symmetric bilinear form on 
$\widetilde{V}$ that has index $n-p$. Then the form ${\tau }_{+}$ on 
${\widehat{V}}_{{\widehat{\sigma }}_{+}}$ is 
nondegenerate and symmetric with index $n-p +1$. The matrix of 
${\tau }_{+}$ with respect to the basis $\widehat{\epsilon} = \{ e_0; 
e_1, \ldots , e_n; e_{n+1} \}$ of ${\widehat{V}}_{{\widehat{\sigma }}_{+}}$ is 
$T_{+}=${\tiny $\begin{pmatrix} 
0 & 0 & 1 \\ 0 & {\widetilde{T}}_{+} & 0 \\ 1 & 0 & 0 \end{pmatrix}$}, where 
${\widetilde{T}}_{+} = I_{n-p,p}=${\tiny $\begin{pmatrix} -I_{n-p} & 0 \\ 0 & I_p 
\end{pmatrix}$} is the matrix of $\widetilde{\tau}$ with respect to the \linebreak 
basis $\widetilde{\epsilon} = \{ e_1, \ldots ,e_n \} $ of $\widetilde{V}$. Let ${\widehat{v}}^{\, 0}$ be a ${\tau }_{+}$ isotropic vector in 
${\widehat{V}}_{{\widehat{\sigma }}_{+}}$, that is, ${\tau }_{+}({\widehat{v}}^{\, 0}, 
{\widehat{v}}^{\, 0}) = 0$. Choose ${\widehat{v}}^{\, 0} = e_{n+1}$. By definition 
${\widehat{v}}^{\, 0}$ is a special vector. A calculation shows that the isotropy group 
\begin{align*}
G(\widehat{V}, {\widehat{\sigma}}_{+}, \tau )_{ {\widehat{v}}^{\, 0}} & = \\
&\hspace{-.75in} = \left\{ \mbox{ \footnotesize $\begin{pmatrix} 1 & 0 & 0 \\ 
\widetilde{d} & \widetilde{A} & 0 \\ -\onehalf {\widetilde{d}}^{\, T} 
{\widetilde{T}}_{+}\widetilde{d} 
& - {\widetilde{d}}^{\, T} \widetilde{T}_{+}\widetilde{A} & 1 \end{pmatrix} 
\in \Oo ({\R }^{n+2}, T_{+}) \setrule \, \widetilde{d} \in {\R }^n \, \, \& \, \, 
\widetilde{A} \in \Oo ({\R}^n, {\widetilde{T}}_{+}) $} \right\} \\
&\hspace{-.75in} = {\Oo ({\R }^{n+2}, T_{+})}_{e_{n+1}} . 
\end{align*}
The mapping
\begin{displaymath}
{G(\widehat{V}, {\widehat{\sigma}}_{+}, \tau )}_{{\widehat{v}}^{\, 0}} \longrightarrow 
G(\widetilde{V}, {\widetilde{\sigma }}_{+}, \widetilde{\tau }) \ltimes \widetilde{V}: 
\mbox{\footnotesize $\begin{pmatrix} 1 & 0 & 0 \\ 
\widetilde{d} & \widetilde{A} & 0 \\ -\onehalf {\widetilde{d}}^{\, T} {\widetilde{T}}_{+}\widetilde{d} 
& - {\widetilde{d}}^{\, T} {\widetilde{T}}_{+}\widetilde{A} & 1 \end{pmatrix} $} 
\longmapsto \mbox{\footnotesize $\begin{pmatrix} 1 & 0 \\ \widetilde{d} & \widetilde{A} \end{pmatrix} $}
\end{displaymath}
is a group isomorphism. The semidirect product 
$G(\widetilde{V}, {\widetilde{\sigma }}_{+}, \widetilde{\tau }) \ltimes \widetilde{V}$ 
is isomorphic to the real affine orthogonal group
$\mathrm{Aff}\Oo (\widetilde{V}, {\widetilde{\sigma }}_{+}, \widetilde{\tau })$ 
of real affine orthogonal mappings of $({\R }^n, {\widetilde{T}}_{+})$ into itself.  \medskip 

Suppose that ${\widetilde{\tau }}_{-} = 
\widetilde{\tau }|{\widetilde{V}}_{{\widetilde{\sigma }}_{+}}$ is nondegenerate 
and alternating. Then $n = \dim \widetilde{V}$ is even, say $2m$. The 
\emph{alternating} bilinear form ${\tau }_{-} = 
\tau |{\widehat{V}}_{{\widehat{\sigma }}_{+}}$, where 
${\widehat{\sigma }}_{+} = {\widehat{\sigma }}^{(n)}_{+}$, with respect to the basis 
${\epsilon }^{\, \raisebox{-4pt}{$\widehat{}$} } = \{ e_0; e_1, \ldots , e_m, f_1, \ldots ,f_m; e_{n+1}\} $\, of 
${\widehat{V}}_{{\widehat{\sigma}}_{+}}$ has 
matrix $T_{-}=${\tiny 
$\begin{pmatrix} 0 & 0 & -1 \\ 0 & {\widetilde{T}}_{-} & 0 \\ 1 & 0 & 0 \end{pmatrix}$}, where ${\widetilde{T}}_{-}=$
{\tiny $\begin{pmatrix} 0 & I_m \\ -I_m & 0 \end{pmatrix} $} is the matrix of 
${\widetilde{\tau }}_{-}$ on ${\widetilde{V}}_{{\widetilde{\sigma }}_{+}}$, which 
has a basis $\{ e_1, \ldots , e_m, f_1, \ldots , f_m \}$. Let ${\widehat{v}}^{\, 0}$ be a nonzero vector in ${\widetilde{V}}_{{\widetilde{\sigma }}_{+}}$. Since ${\tau }_{-}$ is alternating on ${\widehat{V}}_{{\widetilde{\sigma}}_{+}}$, the vector 
${\widehat{v}}^{\, 0}$ is $\tau $ isotropic and hence is a special vector. We may choose ${\widehat{v}}^{\, 0} = e_{n+1}$. Consider the isotropy group 
\begin{align*}
{G(\widehat{V}, {\widehat{\sigma}}_{+}, \tau )}_{{\widehat{v}}^{\, 0}} = 
\left\{ g \in G (\widehat{V}, {\widehat{\sigma}}_{+}, \tau ) \setrule \, 
g{\widehat{\sigma}}_{+} = {\widehat{\sigma}}_{+} g, \, g^{\ast }\tau = \tau  \, \, \& \, \, 
g({\widehat{v}}^{\, 0}) = {\widetilde{v}}^{\, 0} \right\} .
\end{align*}
A calculation gives
\begin{align*}
{G(\widehat{V}, {\widehat{\sigma}}_{+}, \tau )}_{{\widehat{v}}^{\, 0}} = \\
&\hspace{-.75in} = \left\{ \mbox{\footnotesize $\begin{pmatrix} 1 & 0 & 0 \\ \widetilde{d} & \widetilde{A} & 0 \\
0 & -{\widetilde{d}}^{\, T}{\widetilde{T}}_{-}\widetilde{A} & 1 \end{pmatrix} 
\in \Spp ({\R }^{2(m+1)}, T_{-}) \setrule \, \widetilde{d} \in {\R }^{2m} \, \, \& 
\, \, \widetilde{A} \in \Spp ({\R }^{2m}, {\widetilde{T}}_{-}) $} \right\} \\
&\hspace{-.75in} = {\Spp ({\R }^{2(m+1)}, T_{-})}_{e_{n+1}} .
\end{align*}
The mapping 
\begin{displaymath}
{G(\widehat{V}, {\widehat{\sigma}}_{+}, \tau )}_{{\widehat{v}}^{\, 0}} \longrightarrow 
G(\widetilde{V}, {\widetilde{\sigma }}_{+}, \widetilde{\tau }) \ltimes \widetilde{V}: 
\mbox{\footnotesize $\begin{pmatrix} 1 & 0 & 0 \\ 
\widetilde{d} & \widetilde{A} & 0 \\ 0 & - {\widetilde{d}}^{\, T} {\widetilde{T}}_{-}\widetilde{A} & 
1 \end{pmatrix} $} 
\longmapsto \mbox{\footnotesize $\begin{pmatrix} 1 & 0 \\ \widetilde{d} & \widetilde{A} \end{pmatrix} $}
\end{displaymath}
is a group isomorphism. The semidirect product 
$G(\widetilde{V}, {\widetilde{\sigma }}_{+}, \widetilde{\tau }) \ltimes \widetilde{V}$ 
is isomorphic to the real affine symplectic group 
$\mathrm{Aff}\Spp ({\R }^{2m}, {\widetilde{T}}_{-})$ of real affine symplectic mappings of $({\R }^{2m}, {\widetilde{T}}_{-})$ into itself. \medskip

This completes the discussion of the cases where $\widehat{\sigma } = 
{\widehat{\sigma }}_{+}$. \medskip  

\noindent \textsc{Case} $\mathbf{3}$. $\mathrm{Aff}G(\widetilde{V}, {\widetilde{\sigma }}_{-})$. 
Let $V^{\vvee} = \widetilde{V} \times {\C}^2$, where $\widetilde{V}$ is a 
complex vector space of dimension $n =2m$. Set 
\begin{equation}
\begin{array}{l}
{\sigma }^{\vvee}_{-}: V^{\vvee} \rightarrow V^{\vvee}: \\
\hspace{.25in} \big( (\widetilde{u}, \widetilde{w}), (u_{m+1},w_{m+1}) \big)^T \mapsto  
\big( {\widetilde{\sigma }}_{-}(\widetilde{u}, \widetilde{w})^T, 
({\overline{w}}_{m+1}, -{\overline{u}}_{m+1} ) \big) , 
\end{array} 
\label{eq-3one}
\end{equation}
where 
\begin{equation}
{\widetilde{\sigma }}_{-}: \widetilde{V} \rightarrow \widetilde{V}: 
(\widetilde{u}, \widetilde{w})^T \mapsto 
(\overline{\widetilde{w}}, -\overline{\widetilde{u}})^T, 
\label{eq-3two}
\end{equation}
see appendix 1. Then ${\sigma }^{\vvee}_{-}$ is an anti-linear mapping of $V^{\vvee}$ into itself 
such that $({\sigma}^{\vvee})^2_{-} = - {\mathrm{id}}_{V^{\vvee}}$. Let 
${\overline{\sigma }}^{\vvee}_{-}: V^{\vvee} \rightarrow V^{\vvee}$ be the 
complex linear mapping associated to ${\sigma }^{\vvee}_{-}$. In other words, 
${\overline{\sigma }}^{\vvee}_{-}(\alpha v^{\vvee}) = 
\alpha {\sigma }^{\vvee}_{-}(v^{\vvee})$ for every $\alpha \in \C$ and every 
$v^{\vvee} \in V^{\vvee}$. With respect to the basis ${\epsilon }^{\vvee} = 
\{ e_1, \ldots , e_m, f_1, \ldots , f_m; e_{m+1}, f_{m+1} \}$ of $V^{\vvee}$ the 
matrix of ${\overline{\sigma }}^{\vvee}_{-}$ is {\tiny $\left( \begin{array}{cc|cc}
0 & I_m & 0 & 0 \\
-I_m & 0 & 0 & 0 \\ \hline 
0 & 0 & 0 & 1 \\
0 & 0 & -1 & 0 \end{array} \right) $}, using equations (\ref{eq-3one}) and 
(\ref{eq-3two}). Let $\mathcal{A} \in \Gl (V^{\vvee})$. Then $\mathcal{A}$ 
commutes with ${\overline{\sigma }}^{\vvee}_{-}$ if its $2(m+1) \times 2(m+1)$ 
complex matrix with respect to the basis ${\epsilon }^{\vvee}$ is of the form
\begin{displaymath}
\mbox{\footnotesize $\left( \begin{array}{rc|cc} 
A & B & a & b \\
-B & A & -b & a \\ \hline 
\rule{0pt}{10pt} c^T & d^T & e & f \\
-d^T & c^T & -f^T & e \end{array} \right)$,} \quad 
\parbox[t]{2.75in}{where $A$, $B \in \gl ({\C}^m)$; $a$, $b$, $c$, $d \in 
{\C}^m$; \\ $e$, $f \in \C$.}
\end{displaymath}  
Hence $\mathcal{A} \in \Gl (V^{\vvee}, {\sigma }^{\vvee}_{-})$ if and only if 
its $2(m+1) \times 2(m+1)$ complex matrix is of the form 
\begin{displaymath}
\mbox{\footnotesize $\left( \begin{array}{rc|cc} 
A & B & a & b \\
-\overline{B} & \overline{A} & -\overline{b} & \overline{a} \\ \hline 
\rule{0pt}{10pt} c^T & d^T & e & f \\
-{\overline{d}}^{\, T} & {\overline{c}}^{\, T} & -\overline{f} & 
\overline{e} \end{array} \right)$.}
\end{displaymath}
Turn the $2(m+1)$ dimensional complex vector space $V^{\vvee}$ into 
an $m+1$ dimensional quaternionic vector space 
$V^{\vvee}_{{\sigma }^{\vvee}_{-}}$ by defining scalar multiplication $\cdot $ as 
$(\alpha + \beta j) \cdot v^{\vvee} = 
\alpha v^{\vvee} + \beta {\sigma }^{\vvee}_{-}(v^{\vvee})$, 
where $\alpha $, $\beta \in \C$ (and thus $\alpha + \beta j \in \Hh $) and 
$v^{\vvee} \in V^{\vvee}$. The complex linear isomorphism 
\begin{displaymath}
\begin{array}{l}
V^{\vvee}_{{\sigma }^{\vvee}_{-}} \rightarrow V^{\vvee}_{{\sigma }^{\vvee}_{-}}: \\
v^{\vvee} = \big( (\widetilde{u}, \widetilde{w}), (u_{m+1}, w_{m+1}) \big)^T 
\mapsto v^{\vvee} + {\overline{\sigma }}^{\vvee}_{-}(v^{\vvee}) \\
\rule{0pt}{12pt}\hspace{1in} = 
\big( (\widetilde{u} + \widetilde{w}j, \widetilde{w}- \widetilde{u}j), 
(u_{m+1} + w_{m+1}j, w_{m+1} - u_{m+1}j) \big) 
\end{array} 
\end{displaymath}
induces the quaternionic isomorphism 
\begin{displaymath}
\begin{array}{l}
{\rho }^{\vvee}: V^{\vvee}_{{\sigma }^{\vvee}_{-}} \rightarrow 
(V^{\vvee}_{{\sigma }^{\vvee}_{-}})^{\Hh} = {\Hh}^{m+1}: \\
\hspace{.35in} v^{\vvee} = \big( (\widetilde{u}, \widetilde{w}), (u_{m+1}, w_{m+1})^T \big) \mapsto (\widetilde{u}^T + \widetilde{w}^T j, u_{m+1} + w_{m+1}j). 
\end{array} 
\end{displaymath} 
The isomorphism ${\rho }^{\vvee}$ gives rise to the group isomorphism 
\begin{equation}
\begin{array}{l}
\widetilde{\Sigma }: G(V^{\vvee}, {\sigma }^{\vvee}_{-}) \rightarrow 
\Gl \big( (V^{\vvee}_{{\sigma }^{\vvee}_{-}})^{\Hh} \big) : \\
\hspace{.15in} \mathcal{A} = 
 \mbox{\footnotesize $\left( \begin{array}{rc|cc} 
A & B & a & b \\
-\overline{B} & \overline{A} & -\overline{b} & \overline{a} \\ \hline 
\rule{0pt}{10pt} c^T & d^T & e & f \\
-{\overline{d}}^{\, T} & {\overline{c}}^{\, T} & -\overline{f} & 
\overline{e} \end{array} \right)$} \mapsto 
{\rho }^{\vvee} \mathcal{A} ({\rho }^{\vvee})^{-1} = \boldsymbol{\mathcal{A}} = 
\mbox{\footnotesize $\left( \begin{array}{cc} 
{\widetilde{\boldsymbol{A}}}^{\, T} & \widetilde{\boldsymbol{c}} \\ 
{\widetilde{\boldsymbol{a}}}^{\, T} & 
{\widetilde{\boldsymbol{e}}}^{\, T} \end{array} \right) $, }
\end{array}
\label{eq-4threestar}
\end{equation}
where ${\widetilde{\boldsymbol{A}}}^{\, T} = A^T + B^Tj \in \Gl ({\Hh}^m)$; 
${\widetilde{\boldsymbol{a}}}^{\, T} = a^T + b^Tj$, 
$\widetilde{\boldsymbol{c}} = c^T + d^Tj \in {\Hh}^m$; and 
${\widetilde{\boldsymbol{e}}}^{\, T} = e + fj \in \Hh $. 
Note all vectors in ${\Hh }^m$ are \emph{row} vectors and the 
matrices in $\Gl ({\Hh}^m)$ operate on the \emph{right} on ${\Hh }^m$. Moreover, 
multiplication in $\Gl ({\Hh}^{m+2})$ is 
$\boldsymbol{\mathcal{A}}\boldsymbol{\mathcal{B}}$, where 
$\boldsymbol{\mathcal{A}} \in \Gl ({\Hh}^{m+2})$ is applied first followed by 
$\boldsymbol{\mathcal{B}} \in \Gl ({\Hh}^{m+2})$. 
For a special vector $(v^{\vvee})^0 \in V^{\vvee}_{{\sigma }^{\vvee}_{-}}$ let 
${G(V^{\vvee}, {\sigma }^{\vvee})}_{(v^{\vvee})^0}$ be the isotropy group 
consisting of the elements of $G(V^{\vvee}, {\sigma }^{\vvee})$ which leave 
the vector $(v^{\vvee})^0$ fixed. Using the basis ${\epsilon }^{\vvee}$ of 
$V^{\vvee}_{{\sigma }^{\vvee}_{-}}$, choose $(v^{\vvee})^0 = 
e_{m+1}+\mathrm{i} f_{m+1}$. Then 
\begin{displaymath}
{\overline{\sigma }}^{\vvee}_{-}(e_{m+1} 
+ \mathrm{i} f_{m+1}) = -f_{m+1} +\mathrm{i} e_{m+1} = 
\mathrm{i} (e_{m+1} + \mathrm{i} f_{m+1}). 
\end{displaymath}
So $(v^{\vvee})^0$ is a special vector. Moreover, 
\begin{displaymath}
\mathcal{A}((v^{\vvee})^0) = \big( a + \mathrm{i}b ,
-\overline{b} + \mathrm{i}\overline{a} \, \setrule   
e+ \mathrm{i} f , -\overline{f} + \mathrm{i} \overline{e} \big)^T 
= \big( 0 , 0 \, \setrule 1 ,\mathrm{i} \big)^T = (v^{\vvee})^0, 
\end{displaymath}
if and only if $a=0$, $b=0$ and $e=1$, $f=0$. Hence
\begin{align*}
{G(V^{\vvee}, {\sigma }^{\vvee})}_{(v^{\vvee})^0} & = \\
& \hspace{-.9in} = \left\{ \mbox{\footnotesize $\left( \begin{array}{cc|cc} 
A & B & 0 & 0 \\
-\overline{B} & \overline{A} & 0 & 0 \\ \hline 
\rule{0pt}{10pt}c^T & d^T & 1 & 0 \\
-{\overline{d}}^{\, T} & {\overline{c}}^{\, T} & 0 & 1 \end{array} \right) $} 
\in G(V^{\vvee}, {\sigma }^{\vvee}) \setrule \, 
\mbox{\footnotesize $\begin{pmatrix} A & B \\ -\overline{B} & \overline{A} 
\end{pmatrix} $} \in \Gl ({\C}^{2m}) \, \, \& \, \, 
c, \, d \in {\C }^m \right\} . 
\end{align*}
The isomorphism $\widetilde{\Sigma}$ (\ref{eq-4threestar}) restricts to the isomorphism 
\begin{displaymath}
\begin{array}{l}
\Sigma : {G(V^{\vvee}, {\sigma }^{\vvee}_{-})}_{(v^{\vvee})^0} \rightarrow 
\mathrm{Aff}\big( (V^{\vvee}_{{\sigma }^{\vvee}_{-}})^{\Hh} \big) : \\
\hspace{.35in} \mathcal{A} = \mbox{\footnotesize $\left( \begin{array}{rc|cc} 
A & B & 0 & 0 \\
-\overline{B} & \overline{A} & 0 & 0 \\ \hline 
\rule{0pt}{10pt} c^T & d^T & 1 & 0 \\
-{\overline{d}}^{\, T} & {\overline{c}}^{\, T} & 0 & 1 \end{array} \right)$} \mapsto 
{\rho }^{\vvee} \mathcal{A} ({\rho }^{\vvee})^{-1} = \boldsymbol{\mathcal{A}} = 
\mbox{\footnotesize $\left( \begin{array}{cc} 
{\widetilde{\boldsymbol{A}}}^{\, T} & \widetilde{\boldsymbol{c}} \\ 
\widetilde{\boldsymbol{0}} & \widetilde{\boldsymbol{1}} \end{array} \right) $. }
\end{array}
\end{displaymath}

\noindent \textsc{Case} $\mathbf{4}$. $\mathrm{Aff}G(\widetilde{V}, {\widetilde{\sigma }}_{-}, 
\widetilde{\tau })$. Let $\widehat{V} = {\C }^2 \times \widetilde{V} \times {\C }^2$, 
where $\dim \widetilde{V} = 2m$. Define 
\begin{equation}
\begin{array}{l}
{\widehat{\sigma}}_{-}: \widehat{V} \rightarrow \widehat{V}: \\
\hspace{.25in}(z_0,w_0; \widetilde{v}, z_{m+1}, w_{m+1})^T 
\mapsto ({\overline{w}}_0, -{\overline{z}}_0; 
{\widetilde{\sigma }}_{-}(\widetilde{v}); {\overline{w}}_{m+1}, {\overline{z}}_{m+1})^T , 
\end{array}
\label{eq-4one}
\end{equation}
where 
\begin{equation}
{\widetilde{\sigma }}_{-}: \widetilde{V} \rightarrow \widetilde{V}: 
\widetilde{v} = (\widetilde{u}, \widetilde{w})^T \mapsto (\overline{\widetilde{w}}, 
-\overline{\widetilde{u}})^T, 
\label{eq-4two}
\end{equation}
see appendix 1. Then ${\widehat{\sigma }}_{-}$ is an anti-linear mapping 
of $\widehat{V}$ into itself with $({\widehat{\sigma }}_{-})^2 = 
-{\mathrm{id}}_{\widehat{V}}$. Let ${\overline{{\widehat{\sigma }}}}_{-}: 
\widehat{V} \rightarrow \widehat{V}$ be the complex linear mapping 
associated to ${\widehat{\sigma }}_{-}$. With respect to the basis 
$\widehat{\epsilon} = \{ e_0, f_0; e_1, \ldots , e_m, f_1, \ldots , f_m; e_{m+1}, 
f_{m+1} \}$ of $\widehat{V}$ the matrix of ${\overline{\widehat{\sigma }}}_{-}$ is 
\begin{displaymath}
\mbox{\tiny $\left( \begin{array}{cc|cc|cc} 0 & 1 & & & & \\
-1 & 0 & & & & \\ \hline 
& & 0 & I_m & & \\
& & -I_m & 0 & & \\ \hline
& & & & 0 & 1 \\
& & & & -1 & 0 \end{array} \right) $}, 
\end{displaymath}
using equation (\ref{eq-4one}) and (\ref{eq-4two}). Let $\mathcal{A} \in 
\Gl (\widehat{V})$. $\mathcal{A}$ commutes with ${\overline{\widehat{\sigma }}}_{-}$ 
if and only if its $2(m+2) \times 2(m+2)$ complex matrix with respect to the 
basis $\widehat{\epsilon}$ is of the form
\begin{equation}
\mbox{\footnotesize $\left( \begin{array}{c|c|c}
a & b^T & c \\ \hline
\rule{0pt}{10pt}d & \widetilde{A} & e \\ \hline 
\rule{0pt}{10pt}f & g^T & h \end{array} \right) $,}
\label{eq-4three}
\end{equation} 
where $\widetilde{A}=${\tiny $\begin{pmatrix} A & B \\ -B & A \end{pmatrix}$} $\in 
\gl (\widetilde{V})$ with $A$, $B \in \gl ({\C }^m)$; 
$b$, $d$, $e$, $g \in {\mathcal{M}}_{2m,2} = \{ \raisebox{1pt}{ {\tiny $\begin{pmatrix} x_1 & x_2 \\
-x_2 & x_1 \end{pmatrix}$ } } \in M_{2m,2} \setrule \, x_1, \, x_2 \in {\C}^m \} $; 
$a$, $c$, $f$, $h \in {\mathcal{M}}_{2,2} = \{ \raisebox{1pt}{ {\tiny $\begin{pmatrix} 
y_1 & y_2 \\ -y_2 & y_1 \end{pmatrix} $ }} \in \gl ({\C}^2)  
\setrule \, y_1, y_2 \in \C \} $. Here $M_{2m,2}$ is the set of $2m \times 2$ 
complex matrices. Hence $\mathcal{A} \in \Gl (\widehat{V}, {\widehat{\sigma }}_{-})$ if and only if $\mathcal{A}$ is of the form (\ref{eq-4three}), 
where $A=$
{\tiny $\begin{pmatrix} A & B \\ -\overline{B} & \overline{A} \end{pmatrix}$}$\in 
\gl ({\C}^{2m})$ with $A$, $B \in \gl ({\C}^m)$; $b$, $d$, $e$, $g \in 
{\widetilde{\mathcal{M}}}_{2m,2} = $  
$\{ $ \hspace{-5pt} \raisebox{1pt}{\mbox{\tiny $\begin{pmatrix} x_1 & x_2 \\ -{\overline{x}}_2 & {\overline{x}}_1 \end{pmatrix}$} }$\in 
M_{2m,2}(\C) \setrule \, x_1, \, x_2 \in M_{m,2}(\C) \} $; $a$, $c$, $f$, $h \in 
{\widetilde{\mathcal{M}}}_{2,2} = \{ $\raisebox{1pt}{\mbox{\tiny $\begin{pmatrix} 
y_1 & y_2 \\ -{\overline{y}}_2 & {\overline{y}}_1 \end{pmatrix}$}}$\in 
\gl ({\C}^2) \setrule \, y_1,$ $y_2 \in \C \} $. \medskip 

Turn the $2(m+2)$ dimensional complex vector space $\widehat{V}$ into an 
$m+2$ dimensional quaternionic vector space 
${\widehat{V}}_{{\widehat{\sigma }}_{-}}$ by defining scalar multiplication 
$\cdot $ as $(\alpha + \beta j) \cdot \widehat{v} = \alpha \widehat{v} + 
\beta {\widehat{\sigma }}_{-} (\widehat{v})$, 
where $\alpha $, $\beta \in \C$ and $\widehat{v} \in \widehat{V}$. The complex 
linear isomorphism 
\begin{displaymath}
\begin{array}{l}
\widehat{V} \rightarrow {\widehat{V}}^{\Hh}_{{\widehat{\sigma }}_{-}} ={\Hh}^{m+2}: 
\widehat{v} = \big( (u_0, w_0), (\widetilde{u}, \widetilde{w}), 
(u_{m+1}, w_{m+1}) \big)^T \mapsto (\widehat{v} +
{\overline{{\widehat{\sigma }}}}_{-}(\widehat{v}) j)^T \\
\hspace{1in} = \big( (u_0 +w_0j, w_0-u_0j), (\widetilde{u}^{\, T}
+\widetilde{w}^{\, T} j, \widetilde{w}^{\, T} - \widetilde{u}^{\, T} j, \\
\hspace{2.5in}(u_{m+1} +w_{m+1}j, w_{m+1} - u_{m+1}j) \big) 
\end{array}
\end{displaymath}
induces the quaternionic isomorphism 
\begin{equation}
\begin{array}{l}
\widehat{\rho}: {\widehat{V}}_{{\widehat{\sigma}}_{-}} \rightarrow 
{\widehat{V}}^{\Hh}_{{\widehat{\sigma }}_{-}}= {\Hh }^{m+2}: 
\widehat{v} = 
\big( (u_0, w_0), (\widetilde{u}, \widetilde{v}), (u_{m+1}, v_{m+1}) \big)^T  \\ 
\hspace{1.75in}\mapsto (u_0 +w_0j, \widetilde{u}^{\, T} + \widetilde{w}^{\, T} j, 
u_{m+1} +w_{m+1}j ). 
\end{array}
\label{eq-4four}
\end{equation}
This works because 
$(\alpha + \beta j) \cdot \widehat{v} = 
(\alpha + \beta j) \widehat{\rho}(\widehat{v})^T $. \medskip 

Suppose that the bilinear form $\widetilde{\tau }$ on $\widetilde{V}$ is 
\emph{symmetric}. With $\widehat{u} = \big( (u_0,w_0), $ $\widetilde{u}$,  
$(u_{m+1}, w_{m+1}) \big)^T $ and $\widehat{u}' = 
\big( (u'_0,w'_0), \widetilde{u}', (u'_{m+1}, w'_{m+1}) \big)^T \in \widehat{V}$ 
define the bilinear form 
\begin{displaymath}
\tau : \widehat{V} \times \widehat{V}: 
(\widehat{u}, \widehat{u}') \mapsto 
u_0u'_{m+1} + w_0w'_{m+1} + \widetilde{\tau }(\widetilde{u}, \widetilde{u}') 
+ u'_0 u_{m+1} + w'_0w_{m+1} .
\end{displaymath}
With respect to the basis $\widehat{\epsilon}\, ' = 
\{ e'_0, f'_0; e'_1, \ldots , e'_m, f'_1, \ldots , f'_m; e'_{m+1}, f'_{m+1} \}$ of 
$\widehat{V}$ the matrix 
$T_{+}$ of $\tau $ is {\tiny $\left( \begin{array}{r|c|c} 
0 & & I_2 \\ \hline \rule{0pt}{7pt}& {\widetilde{T}}_{+} & \\ \hline 
I_2 & & 0 \end{array} \right) $}, where ${\widetilde{T}}_{+}=${\tiny $\begin{pmatrix}
I_m & 0 \\ 0 & I_m \end{pmatrix}$} is the matrix 
of the bilinear form $\widetilde{\tau}$ on $\widetilde{V}$ with respect to 
the basis $\widetilde{\epsilon}\, '$. Note that for every $\widehat{u}$, 
${\widehat{u}}\, ' \in \widehat{V}$ one has 
$\tau \big( {\widehat{\sigma }}_{-}(\widehat{u}), {\widehat{\sigma }}_{-}(\widehat{u}\, ') 
\big) = \overline{\tau (\widehat{u}, \widehat{u}\, ')}$. On 
${\widehat{V}}_{{\widehat{\sigma}}_{-}} = {\Hh}^{m+2}$ define a quaternion valued form 
\begin{displaymath}
{\tau }_{+}: {\widehat{V}}_{{\widehat{\sigma}}_{-}} \times {\widehat{V}}_{{\widehat{\sigma}}_{-}} 
\rightarrow \Hh : (\widehat{u}, \widehat{v}) 
\mapsto \tau (\widehat{u}, \widehat{v}) + 
\tau (\widehat{u}, {\widehat{\sigma }}_{-}(\widehat{v})) j . 
\end{displaymath}
The form ${\tau }_{+}$ is nondegenerate. For every $\lambda $, $\mu \in \Hh$ one 
has ${\tau }_{+}(\lambda \cdot \widehat{u}, \mu \cdot \widehat{v}) = 
\lambda {\tau }_{+}(\widehat{u}, \widehat{v}) {\mu }^q$ and 
${\tau }_{+}(\widehat{v}, \widehat{u}) = 
{{\tau}_{+}(\widehat{u}, \widehat{v})}^q$, see appendix 1.\footnote{For $\alpha + \beta j \in \Hh$ one has $(\alpha +\beta j)^q = \alpha - \overline{\beta }j$.} There is a basis 
$\epsilon =\{ e_0,f_0; e_1, \ldots , e_m , f_1, \ldots , f_m; e_{m+1}, f_{m+1} \}$ 
of $\widehat{V}$ such that the matrix of ${\tau }_{+}$ with respect to 
$\epsilon $ is {\tiny $\left( \begin{array}{c|c|c} 0 & &I_2 \\ \hline 
& I_{2m} & \\ \hline 
I_2 & & 0 \end{array} \right) $}. On ${\widehat{V}}^{\Hh}_{{\widehat{\sigma }}_{-}}$ 
define a hamiltonian symmetric form ${\boldsymbol{\tau }}_{+}: 
{\widehat{V}}^{\Hh}_{{\widehat{\sigma }}_{-}} \times 
{\widehat{V}}^{\Hh}_{{\widehat{\sigma }}_{-}}  \rightarrow \Hh $ 
by 
\begin{displaymath}
{\boldsymbol{\tau}}_{+}(\widehat{u} + \widehat{v}j, \widehat{w} + \widehat{z}j) 
= (\widehat{u}+ \widehat{v}j)\mbox{\footnotesize $\left( \begin{array}{c|c|c} 
\boldsymbol{0} & & \boldsymbol{1} \\ \hline 
& {\boldsymbol{I}}_m & \\ \hline 
\boldsymbol{1} & & \boldsymbol{0} \end{array} \right) $} 
({\widehat{w}}^T+ {\widehat{z}}^Tj)^q. 
\end{displaymath}
The matrix of ${\boldsymbol{\tau}}_{+}$ with respect to the basis 
\begin{displaymath}
\boldsymbol{\eta} = 
\{ \ttfrac{1}{\sqrt{2}}(e_0 +f_0j); \ttfrac{1}{\sqrt{2}}(e_1 + f_1j), \ldots , 
\ttfrac{1}{\sqrt{2}}(e_m+f_mj); \ttfrac{1}{\sqrt{2}}(e_{m+1}+f_{m+1}j) \} 
\end{displaymath}
of ${\widehat{V}}^{\Hh}_{{\widehat{\sigma }}_{-}}$ is 
{\tiny $ \left( \begin{array}{c|c|c} \boldsymbol{0} & & \boldsymbol{1} \\ \hline  
& {\boldsymbol{I}}_m & \\ \hline 
\boldsymbol{1} & & \boldsymbol{0} \end{array} \right) $}. Observe that 
the isomorphism $\widehat{\rho }$ (\ref{eq-4four}) pulls back 
the hamiltonian symmetric form ${\boldsymbol{\tau }}_{+}$ on 
${\widehat{V}}^{\Hh}_{{\widehat{\sigma }}_{-}} = {\Hh}^{m+2}$ to 
the symmetric form ${\tau }_{+}$ on ${\widehat{V}}_{{\widehat{\sigma }}_{-}}$. 
\medskip 

The matrix $\mathcal{A} \in G(\widehat{V}, {\widehat{\sigma }}_{-})$ (\ref{eq-4three}) fixes the vector $(\widehat{v})^0 = e_{m+1} +\mathrm{i}f_{m+1} \in 
{\widehat{V}}_{{\widehat{\sigma }}_{-}}$, which special 
since ${\widehat{\sigma }}_{-}((\widehat{v})^0) = \mathrm{i} (\widehat{v})^0$ 
and $\tau \big( (\widehat{v})^0, (\widehat{v})^0 \big) =0$, if and only if 
$\mathcal{A} =${\tiny 
$\begin{pmatrix} a & b^T & 0 \\ d & \widetilde{A} & 0 \\
f & g^T & I_2 \end{pmatrix}$}. Here $a$, $f$; $b$, $d$, $g$; and 
$\widetilde{A}$ satisfy the conditions following equation (\ref{eq-4three}). A 
calculation shows that $\mathcal{A} \in 
{G(\widehat{V}, {\widehat{\sigma }}_{-})}_{(\widehat{v})^0}$ preserves the 
symmetric bilinear form ${\tau }_{+}$ on $\widehat{V}$ if and only if 
\begin{displaymath}
\mathcal{A} = \mbox{\footnotesize $\begin{pmatrix} 
I_2 & 0 & 0 \\ \widetilde{d} & \widetilde{A} & 0 \\ 
-\onehalf {\widetilde{d}}^{\, T}{\widetilde{T}}_{+}\widetilde{d} & 
-{\widetilde{d}}^{\, T} {\widetilde{T}}_{+} \widetilde{A} & I_2 
\end{pmatrix} $}, 
\end{displaymath}
where $\widetilde{d}=${\tiny $\begin{pmatrix} {\widetilde{d}}_1 & 
{\widetilde{d}}_2 \\ -{\overline{\widetilde{d}}}_2 & {\overline{\widetilde{d}}}_1 
\end{pmatrix}$} with ${\widetilde{d}}_1$, ${\widetilde{d}}_2 \in {\C}^m$ and 
$\widetilde{A} \in G(\widetilde{V}, {\widetilde{\sigma }}_{-}, {\widetilde{\tau}}_{+})$, 
that is, $\widetilde{A} \in \Gl (\widetilde{V})$ such that $\widetilde{A} 
{\widetilde{\sigma }}_{-} = {\widetilde{\sigma }}_{-}\widetilde{A}$ and 
${\widetilde{T}}_{+} = {\widetilde{A}}^{\, T} {\widetilde{T}}_{+} \widetilde{A}$. 
The isomorphism $\widehat{\rho}$ (\ref{eq-4four}) gives rise to the group isomorphism 
\begin{equation}
\begin{array}{l}
\widetilde{\Sigma }: G(\widehat{V}, {\widehat{\sigma }}_{-}, {\tau }_{+}) 
\rightarrow G( {\widehat{V}}^{\Hh}_{{\widehat{\sigma }}_{-}}, 
{\boldsymbol{\tau }}_{+}): \\
\rule{.5in}{0in}\mathcal{A} = \mbox{\footnotesize $\left( \begin{array}{c|c|c}
a & b^T & c \\ \hline
\rule{0pt}{10pt}d & \widetilde{A} & e \\ \hline 
\rule{0pt}{10pt}f & g^T & h \end{array} \right) $} \mapsto 
\widehat{\rho }\mathcal{A} {\widehat{\rho }}^{-1} = \boldsymbol{\mathcal{A}} = 
\mbox{\footnotesize $\left( \begin{array}{c|c|c}
{\boldsymbol{a}}^T & {\boldsymbol{d}}^T & {\boldsymbol{f}}^T\\ \hline
\rule{0pt}{12pt}\boldsymbol{b} & {\boldsymbol{\widetilde{A}}}^T & 
\boldsymbol{g} \\ \hline 
\rule{0pt}{10pt}{\boldsymbol{c}}^T& {\boldsymbol{e}}^T & {\boldsymbol{h}}^T \end{array} \right) $, } \end{array}
\label{eq-4five}
\end{equation} 
where ${\boldsymbol{a}}^T = a_1+a_2j$, ${\boldsymbol{c}}^T = c_1+c_2j$, 
${\boldsymbol{f}}^T = f_1 +f_2j$, and $\boldsymbol{h} = h_1+h_2j \in \Hh$; 
$\boldsymbol{b} = b^T_1 +b^T_2j$, ${\boldsymbol{d}}^T = d^T_1 +d^T_2j$, 
${\boldsymbol{e}}^T = e^T_1 +e^T_2 j$, and $\boldsymbol{g} = g^T_1 +g^T_2 
\in {\Hh}^m$; $\boldsymbol{\widetilde{A}}^T = A^T + B^Tj$ with 
$A$, $B \in \gl ({\C}^m)$. Moreover, 
$(\widetilde{\Sigma })^{\ast} \boldsymbol{\tau }_{+} = {\tau }_{+}$. The map 
$\widetilde{\Sigma }$ (\ref{eq-4five}) restricts to the group 
isomorphism 
\begin{displaymath}
\begin{array}{l}
\Sigma : {G(\widehat{V}, {\widehat{\sigma }}_{-}, {\tau }_{+})}_{(\widehat{v})^0} \rightarrow 
\mathrm{Aff}G({\Hh}^m, {\widetilde{\boldsymbol{\tau }}}_{+}): \\
\rule{.75in}{0in}\mbox{\footnotesize $\begin{pmatrix} 
I_2 & 0 & 0 \\ \widetilde{d} & \widetilde{A} & 0 \\ 
-\onehalf {\widetilde{d}}^{\, T}{\widetilde{T}}_{+}\widetilde{d} & 
-{\widetilde{d}}^{\, T} {\widetilde{T}}_{+} \widetilde{A} & I_2 
\end{pmatrix} $} \mapsto \mbox{\footnotesize $\begin{pmatrix} 
{\widetilde{\boldsymbol{A}}}^{\, T} & {\widetilde{\boldsymbol{d}}}^{\, T} \\ 
\boldsymbol{0} & \boldsymbol{1} 
\end{pmatrix}$,}
\end{array} 
\end{displaymath}
where ${\widetilde{\boldsymbol{A}}}^T = A^T + B^Tj \in 
G({\Hh}^m, {\boldsymbol{\tau }}_{+})$ 
and $\widetilde{\boldsymbol{d}}^{\, T} = {\widetilde{d}}^{\, T}_1 
+{\widetilde{d}}^{\, T}_2 \in \Hh $.  \medskip 

Suppose that $\widetilde{\tau }$ is a nondegnerate \emph{alternating} 
bilinear form on $\widetilde{V}$. With $\widetilde{u} = (\widetilde{z}, \widetilde{w})$, 
${\widetilde{u}}^{\, \prime } = ({\widetilde{z}}^{\, \prime }, {\widetilde{w}}^{\, \prime }) 
\in \widetilde{V}$ let
\begin{displaymath}
\begin{array}{l}
\tau : \widehat{V} \times \widehat{V} \rightarrow \C :
 \big( (\widehat{z}, \widehat{w}), ({\widehat{z}}^{\, \prime }, 
{\widehat{w}}^{\, \prime}) \big) \longrightarrow \\
\hspace{1in}z_0w'_{m+1} -w_0z'_{m+1} + \widetilde{\tau }(\widetilde{u}, {\widetilde{u}}^{\, \prime })+w'_0z_{m+1}- z'_0w_{m+1} . 
\end{array}
\end{displaymath}
Then $\tau $ is a nondegenerate alternating bilinear form on $\widehat{V}$. With 
respect to the basis $\{ e_0, f_0; e_1 \ldots , e_m, f_1, \ldots , f_m; e_{m+1}, f_{m+1} \} $ of $\widehat{V}$ the matrix of $\tau $ is $T_{-}=${\tiny $\left( \begin{array}{c|c|c} 
0 & & J_2 \\ \hline & \rule{0pt}{7pt} {\widetilde{T}}_{-} & \\ \hline J_2 & & 0 \end{array} \right) $}, where ${\widetilde{T}}_{-}= J_{2m-p}=${\tiny $\begin{pmatrix} 
0 & I_{m-p/2} \\ -I_{m-p/2} & 0 \end{pmatrix}$} is the matrix of the nondegenerate alternating form $\widetilde{\tau}$ on $\widetilde{V}$ with respect to the basis 
$\{ e_1, \ldots , e_m; f_1, \ldots , f_m \} $ and $J_2=$\raisebox{1pt}{{\tiny $\begin{pmatrix} 0 & 1 \\ -1 & 0 \end{pmatrix}$}}. Note that 
$\tau \big( {\widehat{\sigma }}_{-}(\widehat{u}), {\widehat{\sigma }}_{-}(\widehat{v}) \big) = \overline{\tau (\widehat{u}, \widehat{v})}$. 
Let ${\tau }_{-} = \tau |{\widehat{V}}_{{\widehat{\sigma }}_{-}}$. Then ${\tau }_{-}$ 
is a nondegenerate alternating bilinear form on 
${\widehat{V}}_{{\widehat{\sigma }}_{-}}$. On 
${\widehat{V}}^{\Hh}_{{\widehat{\sigma }}_{-}} = {\Hh}^{m+2}$ for each 
$0 \le p/2 \le m$ define a hamiltonian alternating quaternion valued form 
${\boldsymbol{\tau }}_{-}: 
{\Hh}^{m+2} \times {\Hh }^{m+2} \rightarrow \Hh $ by 
\begin{displaymath}
{\boldsymbol{\tau }}_{-}(\widehat{u} + \widehat{v}j, \widehat{w} + \widehat{z}j) = 
(\widehat{u} + \widehat{v}j)\mbox{\footnotesize $\left( \begin{array}{c|c|c} 
0 & & j \\ \hline 
& j I_{m-p/2, p/2} & \\ \hline 
j & & 0 \end{array} \right)$} ((\widehat{w})^T + (\widehat{z})^Tj)^q, 
\end{displaymath}
where $I_{m-p/2, p/2} =${\tiny $\begin{pmatrix} -I_{m-p/2} & 0 \\ 0 & I_p \end{pmatrix} $}. The matrix of ${\boldsymbol{\tau }}_{-}$ with respect to the 
basis 
\begin{displaymath}
\boldsymbol{\eta } = \{ \ttfrac{1}{\sqrt{2}}(e_0 +f_0j; 
\ttfrac{1}{\sqrt{2}}(e_1 +f_1j), \ldots ,  \ttfrac{1}{\sqrt{2}}(e_m +f_mj) ; 
\ttfrac{1}{\sqrt{2}}(e_{m+1} +f_{m+1}j) \} 
\end{displaymath}
of ${\widehat{V}}^{\Hh}_{{\widehat{\sigma }}_{-}}$ is 
{\tiny $\left( \begin{array}{c|c|c} 
0 & & j \\ \hline & I_{m-p/2,p/2} j & \\ \hline j & & 0 \end{array} \right)$}. 
Observe that the map $\widehat{\rho}$ (\ref{eq-4four}) pulls back 
the hamiltonian alternating form ${\boldsymbol{\tau}}_{-}$ on 
${\widehat{V}}^{\Hh}_{{\widehat{\sigma }}_{-}}$ to the alternating form 
${\tau }_{-}$ on $\widehat{V}$. \medskip
 
Let $(\widehat{v})^0 = e_{m+1}+\mathrm{i} f_{m+1} \in 
{\widehat{V}}_{{\widehat{\sigma }}_{-}}$. The vector $(\widehat{v})^0$ is 
special since ${\widehat{\sigma}}_{-}\big( (\widehat{v})^0 \big) = 
\mathrm{i} (\widehat{v})^0$ and $\tau \big( (\widehat{v})^0, (\widehat{v})^0 \big) 
=0$, because $\tau $ is alternating. An element $\widehat{A}$ of the isotropy group 
${G(\widehat{V}, {\widehat{\sigma }}_{-})}_{(\widehat{v})^0}$ preserves 
the alternating form ${\tau }_{-}$ on ${\widehat{V}}_{{\widehat{\sigma }}_{-}}$ if 
and only if $\widehat{A} =${\tiny $\begin{pmatrix} I_2 & 0 & 0 \\
d & \widetilde{A} & 0 \\ \onehalf J_2d^T {\widetilde{T}}_{-}d & 
- J_2d^T{\widetilde{T}}_{-} & I_2 \end{pmatrix} $}, where $d=${\tiny 
$\begin{pmatrix} d_1 & d_2 \\ -{\overline{d}}_2 & {\overline{d}}_1 \end{pmatrix}$} 
with $d_1$, $d_2 \in {\C}^m$; and $\widetilde{A} \in G(\widetilde{V}, 
{\widetilde{\sigma}}_{-}, \widetilde{\tau })$, that is, $\widetilde{A} \in 
\Gl (\widetilde{V})$, ${\widetilde{\sigma}}_{-}\widetilde{A} = 
\widetilde{A}{\widetilde{\sigma }}_{-}$, and ${\widetilde{T}}_{-} = 
{\widetilde{A}}^{\, T} {\widetilde{T}}_{-}\widetilde{A}$. The isomorphism 
$\widehat{\rho}$ (\ref{eq-4four}) gives rise to the group isomorphism 
\begin{equation}
\begin{array}{l}
\widetilde{\Sigma }: G(\widehat{V}, {\widehat{\sigma }}_{-}, {\tau }_{-}) 
\rightarrow G( {\widehat{V}}^{\Hh}_{{\widehat{\sigma }}_{-}}, 
{\boldsymbol{\tau }}_{-}): \\
\rule{.5in}{0in}\mathcal{A} = 
\mbox{\footnotesize $\left( \begin{array}{c|c|c}
a & b^T & c \\ \hline
\rule{0pt}{10pt}d & \widetilde{A} & e \\ \hline 
\rule{0pt}{10pt}f & g^T & h \end{array} \right) $} \mapsto 
\widehat{\rho }\mathcal{A} {\widehat{\rho }}^{-1} = \boldsymbol{\mathcal{A}} = 
\mbox{\footnotesize $\left( \begin{array}{c|c|c}
\boldsymbol{a}^T & \boldsymbol{d}^T & \boldsymbol{f}^T \\ \hline
\rule{0pt}{12pt}\boldsymbol{b} & {\boldsymbol{\widetilde{A}}}^T & 
\boldsymbol{g} \\ \hline 
\rule{0pt}{10pt}\boldsymbol{c}^T & \boldsymbol{e}^T & \boldsymbol{h}^T 
\end{array} \right) $, } \end{array}
\label{eq-4six}
\end{equation} 
where $\boldsymbol{a} = a_1+a_2j$, $\boldsymbol{c} = c_1+c_2j$, 
$\boldsymbol{f} = f_1 +f_2j$, and $\boldsymbol{h} = h_1+h_2j \in \Hh$; 
$\boldsymbol{b}^T = b^T_1 +b^T_2j$, $\boldsymbol{d}^T = d^T_1 +d^T_2j$, 
$\boldsymbol{e}^T = e^T_1 +e^T_2 j$, and $\boldsymbol{g}^T = g^T_1 +g^T_2 
\in {\Hh}^m$; $\boldsymbol{\widetilde{A}} = A^T + B^Tj$ with 
$A$, $B \in \gl ({\C}^m)$. Moreover, 
$(\widehat{\rho})^{\ast} {\boldsymbol{\tau }}_{-} = {\tau }_{-}$. The map 
$\widetilde{\Sigma }$ (\ref{eq-4five}) restricts to the group 
isomorphism 
\begin{displaymath}
\begin{array}{l}
\Sigma : {G(\widehat{V}, {\widehat{\sigma }}_{-}, {\tau }_{-})}_{(\widehat{v})^0} \rightarrow 
\mathrm{Aff}G({\Hh}^m, {\widetilde{\boldsymbol{\tau }}}_{-}):\\
\rule{.75in}{0in} \mathcal{A} = \mbox{\footnotesize $\begin{pmatrix} 
I_2 & 0 & 0 \\ \widetilde{d} & \widetilde{A} & 0 \\ 
-\onehalf J_2{\widetilde{d}}^{\, T}{\widetilde{T}}_{-}\widetilde{d} & 
-J_2{\widetilde{d}}^{\, T} {\widetilde{T}}_{-} \widetilde{A} & I_2 
\end{pmatrix} $} \mapsto \boldsymbol{\mathcal{A}}= 
\mbox{\footnotesize $\begin{pmatrix} 
{\widetilde{\boldsymbol{A}}}^{\, T} & {\widetilde{\boldsymbol{d}}}{\, ^T} \\ \boldsymbol{0} & \boldsymbol{1} 
\end{pmatrix}$,} 
\end{array}
\end{displaymath}
where $\widetilde{\boldsymbol{A}}^T = A^T + B^Tj \in 
G({\Hh}^m, {\boldsymbol{\tau }}_{-})$ 
and ${\widetilde{\boldsymbol{d}}}^{\, T} = {\widetilde{d}}^{\, T}_1 +{\widetilde{d}}^{\, T}_2 \in {\Hh }^m$.  \medskip 

This completes case 4 and the proof of theorem 4.1. \hfill $\square $ 

 \bigskip

\noindent {\Large \textbf{Appendix 1. Classification of anti-linear mappings}} \bigskip 
%%%%%%%%%%%%%%%%%%%

Two anti-linear mappings $\widetilde{\sigma }$ and ${\widetilde{\sigma }}'$ of the complex vector space $\widetilde{V}$ into itself are \emph{equivalent} if and only if ${\widetilde{\sigma }}' = \alpha k^{-1}\widetilde{\sigma }k$ for some $k \in G$ and 
some $\alpha \in \C \setminus \{ 0 \}$. Note that $\widetilde{\sigma }$ and 
$-\widetilde{\sigma }$ are equivalent. We may suppose that 
${\widetilde{\sigma }}^2 = \pm {\mathrm{id}}_{\widetilde{V}}$. Equivalence classes of anti-linear mappings for $G = G(\widetilde{V}, \widetilde{\sigma}, \widetilde{\tau })$ 
are determined in \cite[app 1, p.357]{burgoyne-cushman}. We 
list explicit representatives below. \medskip %

\noindent \textsc{Case} $\mathbf{1}$. $\widetilde{\sigma } = 
{\widetilde{\sigma }}_{+}$. \medskip

Suppose that  $\widetilde{V} = {\C }^n$, 
where $z = (z_1, \ldots , z_n)^T \in {\C}^n$. Let 
\begin{displaymath}
\widetilde{\sigma } : {\C}^n \rightarrow {\C}^n: z \mapsto \overline{z} = 
({\overline{z}}_1, \ldots , {\overline{z}}_n)^T. 
\end{displaymath}
Set ${\widetilde{\sigma }}_{+} = 
\widetilde{\sigma }$. Since
\begin{displaymath}
{\widetilde{\sigma }}_{+}(\alpha z + \beta z') = 
\overline{\alpha z + \beta z'} = \overline{\alpha }\, \overline{z} + 
\overline{\beta } \, \overline{z'} = \overline{\alpha }\, {\widetilde{\sigma }}_{+}(z) + 
\overline{\beta }\, {\widetilde{\sigma }}_{+}(z'), 
\end{displaymath}
the mapping ${\widetilde{\sigma }}_{+}$ is anti-linear. Also 
${\widetilde{\sigma }}^{\, 2}_{+} = {\mathrm{id}}_{{\C}^n}$, since 
${\widetilde{\sigma }}^{\, 2}_{+}(z) = 
{\widetilde{\sigma }}_{+}(\overline{z}) = z$. This handles the case 
when $G = \Gl (\widetilde{V}, {\widetilde{\sigma }}_{+})$. 
\par Suppose that on ${\C}^n$ we have 
a nondegenerate complex valued bilinear form 
$\widetilde{\tau }: {\C }^n \times {\C }^n \rightarrow \C$. \medskip 

If $\widetilde{\tau }$ is \emph{symmetric}, we may assume that its matrix with 
respect to the standard basis of ${\C}^n$ is $I_n$, that is, 
$\widetilde{\tau }(z,w) = w^T I_n z$. Let 
\begin{displaymath}
{\widetilde{V}}_{{\widetilde{\sigma }}_{+}} = 
\{ z \in {\C }^n \setrule \, {\widetilde{\sigma}}_{+}(z) = z \} = 
\{ x \in {\R }^n \setrule \, z = x + \mathrm{i}y \} = {\R }^n
\end{displaymath} 
and set ${\widetilde{\tau}}_{+} = 
\widetilde{\tau }|{\widetilde{V}}_{{\widetilde{\sigma }}_{+}}: 
{\R }^n \times {\R }^n \rightarrow \R$. Then ${\widetilde{\tau}}_{+}$ is a 
nondegenerate real valued symmetric bilinear form on ${\R }^n$. Replacing 
${\widetilde{\sigma }}_{+}$ by 
$-{\widetilde{\sigma }}_{+}$, if necessary, we may assume that 
${\widetilde{\tau}}_{+}$ has signature $(n-p,p)$ for 
$0 \le p \le \left[ \frac{n}{2} \right] $. We may suppose that the matrix of 
${\widetilde{\tau }}_{+}$ with respect to the standard basis of ${\R }^n$ is 
$I_{n-p,p}=$\raisebox{1pt}{{\tiny $\begin{pmatrix} -I_{n-p} & 0 \\ 0 & I_p 
\end{pmatrix}$}}. Then 
${\widetilde{\tau }}_{+}(x,y) = y^T I_{n-p,p}x$. Let 
\begin{displaymath}
{\widetilde{\sigma }}^{\, (p)}_{+}: 
{\C }^n \rightarrow {\C}^n : z \mapsto I_{n-p,p}\overline{z}. 
\end{displaymath}
Then %
\begin{align*}
\widetilde{\tau }({\widetilde{\sigma }}^{\, (p)}_{+}(z), 
{\widetilde{\sigma }}^{\, (p)}_{+}(w)) & = (I_{n-p,p}\overline{w})^TI_n (I_{n-p,p}\overline{z}) \\
& = {\overline{w}}^{\, T} \big( (I_{n-p,p})^T I_n I_{n-p,p} \big)\overline{z} 
 = {\overline{w}}^T I_n \overline{z} = \overline{\widetilde{\tau }(z, w)}. 
\end{align*}  
This handles the case when $G = \Oo (\widetilde{V}, {\widetilde{\sigma }}^{(p)}_{+}, 
\widetilde{\tau })$. 
\par If $\widetilde{\tau }$ is \emph{alternating}, then $n$ is even, say $2m$, since 
$\widetilde{\tau }$ is nondegenerate. We may assume that the matrix of the nondegenerate real valued alternating form ${\widetilde{\tau }}_{+}: 
{\R }^m \times {\R }^m \rightarrow \R $ with respect to the standard basis of 
${\R }^n = {\R}^m \times {\R }^m$ is $J_{2m} =$
\raisebox{1pt}{ {\tiny $\begin{pmatrix} 0 & I_m \\ -I_m & 0 \end{pmatrix} $} }, 
that is, $\widetilde{\tau }_{+}\left( \mbox{\tiny $\begin{pmatrix} x \\ y 
\end{pmatrix}$}, \mbox{\tiny $\begin{pmatrix} x' \\ y' \end{pmatrix}$} \right) = 
(x', y')^T J_{2m}${\footnotesize $\begin{pmatrix} x \\ y \end{pmatrix}$}. 
Now 
\begin{align*}
\widetilde{\tau }({\widetilde{\sigma }}_{+}(z,w) , {\widetilde{\sigma }}_{+}(z' , w' ) 
& = \big( ({\overline{z}}\, ' )^T, ({\overline{w}}\, ')^T \big)J_{2m} 
\mbox{\tiny $\begin{pmatrix}\overline{z} \\ \overline{w} \end{pmatrix} $} = 
\overline{\widetilde{\tau }((z,w), (z',w'))}. 
\end{align*}
This handles the case when $G = \Spp (\widetilde{V}, {\widetilde{\sigma}}_{+}, 
\widetilde{\tau })$. \medskip   

Suppose that ${\widetilde{\sigma }}^2 = -{\mathrm{id}}_{\widetilde{V}}$.  Let 
$\mathbb{H} = \{ \alpha + \beta j \setrule \, \alpha , \beta \in \C; \, 
j^2 =-1; \, \alpha j = j \overline{\alpha }\} $ be the quaternions with anti-involution 
$(\alpha + \beta j)^q = \alpha  - \overline{\beta }j$. Note that if $x$ and 
$y \in \mathbb{H}$, then $(xy)^q = y^q x^q$. Let $\widetilde{V} = {\C}^n$ with $n =2m$. Consider the anti-linear mapping 
\begin{displaymath}
\widetilde{\sigma }: {\C}^n = {\C}^m \times {\C}^m \rightarrow {\C }^n: 
(z,w)^T \mapsto (-\overline{w}, \overline{z})^T .
\end{displaymath}
Then ${\widetilde{\sigma }}^{\, 2} = -{\mathrm{id}}_{{\C}^n}$. Set 
${\widetilde{\sigma }}_{-} = \widetilde{\sigma }$. Turn ${\C}^n$ into a vector space 
${\widetilde{V}}_{{\widetilde{\sigma }}_{-}}$ over the quaternions by 
defining scalar multiplication $\cdot $ as $(\alpha + \beta j)\cdot (z,w) = 
\alpha (z,w) + \beta {\widetilde{\sigma }}_{-}(z,w)$. The map 
$\rho : {\widetilde{V}}_{{\widetilde{\sigma }}_{-}} \rightarrow {\mathbb{H}}^m: 
(z,w)^T \mapsto z^T + w^T j $ is an isomorphism of quaternionic vector spaces, since 
\begin{align*}
\rho \big( (\alpha + \beta j)\cdot (z,w)^T \big) & = 
\rho \big( \alpha (z,w)^T + \beta {\widetilde{\sigma }}_{-}((z,w)^T) \big) \\
& = (\alpha z^T - \beta \overline{w}^T) + (\alpha w^T +\beta \overline{z}^T)j 
= (\alpha + \beta j)(z^T +w^T j). 
\end{align*}
This takes case of the case when $G = 
\Gl (\widetilde{V}, {\widetilde{\sigma }}_{-})$. \medskip

Suppose that $\widetilde{\tau }$ is a nondegenerate complex valued 
bilinear form on ${\C}^n$. On ${\widetilde{V}}_{{\widetilde{\sigma }}_{-}}$ 
define the quaternion valued bilinear form
\begin{displaymath}
{\widetilde{\tau }}_{-}: {\widetilde{V}}_{{\widetilde{\sigma }}_{-}} \times 
{\widetilde{V}}_{{\widetilde{\sigma }}_{-}} \rightarrow \mathbb{H}:
(\widetilde{u}, \widetilde{v}) \mapsto \widetilde{\tau }(\widetilde{u}, \widetilde{v}) 
+ \widetilde{\tau}\big( \widetilde{u}, {\widetilde{\sigma }}_{-}(\widetilde{v}) \big) j 
\end{displaymath}
Using the fact that $\widetilde{\tau}\big( {\widetilde{\sigma }}_{-}(\widetilde{u}), 
{\widetilde{\sigma }}_{-}(\widetilde{v}) \big) = 
\overline{\widetilde{\tau}(\widetilde{u}, \widetilde{v})}$, 
which implies 
$\widetilde{\tau}\big( {\widetilde{\sigma }}_{-}(\widetilde{u}), 
\widetilde{v}) \big) = - \overline{\widetilde{\tau}\big( \widetilde{u}), 
{\widetilde{\sigma }}_{-}(\widetilde{v}) \big)} $. 
A straightforward calculations shows that for every $\lambda \in \mathbb{H}$ 
\begin{subequations}
\begin{equation}
{\widetilde{\tau }}_{-}\big( \lambda \cdot \widetilde{u}, \widetilde{v} \big) = 
\lambda {\widetilde{\tau }}_{-}(\widetilde{u}, \widetilde{v}) 
\label{eq-app1onea}
\end{equation}
and for every $\mu \in \mathbb{H}$ 
\begin{equation}
{\widetilde{\tau }}_{-}\big(\widetilde{u}, \mu \cdot \widetilde{v} \big) = 
{\widetilde{\tau }}_{-}(\widetilde{u}, \widetilde{v}) {\mu }^q. 
\label{eq-app1oneb}
\end{equation}
\end{subequations}
 
Suppose that $\widetilde{\tau }$ is \emph{symmetric}. We may assume that the 
matrix of $\widetilde{\tau}$ with respect to the standard basis of 
${\C }^m \times {\C}^m = {\C }^n$ is \raisebox{2pt}{{\tiny $\begin{pmatrix} I_m & 0 \\
0 & I_m \end{pmatrix}$}}, that is, $\widetilde{\tau }\big( \widetilde{u},$ $
\widetilde{v} \big) $ 
$ = {\widetilde{v}}^{\, T} ${\tiny $\begin{pmatrix} I_m & 0 \\ 0 & I_m 
\end{pmatrix}$}$\widetilde{u}$. Now 
\begin{align*}
\widetilde{\tau } \big( {\widetilde{\sigma }}_{-}(\widetilde{u}, 
{\widetilde{\sigma }}_{-}(\widetilde{v}) \big) & = 
\mbox{\footnotesize $\left( \begin{pmatrix} 0 & -I_m \\ I_m & 0 \end{pmatrix} \right. $} \left. \! \! {\widetilde{v}}\, \right)^{T}
\mbox{\footnotesize $\begin{pmatrix} I_m & 0 \\ 0 & I_m \end{pmatrix}$} 
\mbox{\footnotesize $ \begin{pmatrix} I_m & 0 \\ 0 & I_m \end{pmatrix}$}
\widetilde{u} \\
&= {\widetilde{v}}^{\, T}\mbox{\footnotesize $\begin{pmatrix} I_m & 0 \\ 0 & I_m 
\end{pmatrix}$} \widetilde{u} = 
\overline{\widetilde{\tau }(\widetilde{u}, \widetilde{v})} , 
\end{align*}
which implies 
\begin{align*}
\widetilde{\tau }\big( \widetilde{u}, {\widetilde{\sigma }}_{-}(\widetilde{v}) \big) & = 
- \widetilde{\tau }\big( {\widetilde{\sigma }}_{-}( {\widetilde{\sigma }}_{-}(\widetilde{u}) ), {\widetilde{\sigma }}_{-}(\widetilde{v}) \big) 
= -\overline{\widetilde{\tau }\big( {\widetilde{\sigma }}_{-}(\widetilde{u}, 
\widetilde{v}) \big)} . 
\end{align*}
So
\begin{align*}
{{\widetilde{\tau }}_{-}( \widetilde{u}, \widetilde{v})}^q & = 
\big( \widetilde{\tau }\big( \widetilde{u}), \widetilde{v} \big) + 
\widetilde{\tau } \big( \widetilde{u}, {\widetilde{\sigma }}_{-}(\widetilde{v}) \big) j \big)^q \\
& = \widetilde{\tau } ( \widetilde{u}, \widetilde{v} ) 
- \overline{\widetilde{\tau}\big( \widetilde{u}, {\widetilde{\sigma }}_{-}(\widetilde{v})\big)} j 
 = \widetilde{\tau } (\widetilde{u}, \widetilde{v} ) + 
\widetilde{\tau }\big( {\widetilde{\sigma}}_{-}(\widetilde{u}), \widetilde{v} \big) j \\
& = \widetilde{\tau }( \widetilde{v}, \widetilde{u}) + 
\widetilde{\tau }\big( (\widetilde{v}, {\widetilde{\sigma}}_{-}(\widetilde{u}) \big) j 
= {\widetilde{\tau }}_{-}(\widetilde{v}, \widetilde{u}) ,  
\end{align*}
that is, ${\widetilde{\tau }}_{-}$ is \emph{hamiltonian symmetric}. Let 
$\epsilon = {\{ \frac{1}{\sqrt{2}}(e_{\ell }, f_{\ell }) \}}^m_{\ell =1}$ be a 
$\widetilde{\tau }$ orthogonal basis of 
${\widetilde{V}}_{{\widetilde{\sigma }}^{(p)}_{-}}$. Then 
\begin{align*}
{\widetilde{\tau }}_{-}\big(  \ttfrac{1}{\sqrt{2}}(e_{\ell}, f_{\ell}),  
\ttfrac{1}{\sqrt{2}}(e_k, f_k) ) & = 
\onehalf \widetilde{\tau }\big( (e_{\ell}, f_{\ell }), (e_k, f_k) \big) + 
\widetilde{\tau }\big( (e_{\ell }, f_{\ell }), {\widetilde{\sigma }}_{-}(e_k, f_k) \big )j \\
& = \onehalf \widetilde{\tau }\big( (e_{\ell}, f_{\ell }), (e_k, f_k) \big) + 
\widetilde{\tau }\big( (e_{\ell }, f_{\ell }), (-f_k, e_k) \big )j \\
& = {\delta }_{\ell k}. 
\end{align*}
In other words, with respect to the basis $\epsilon $ of ${\mathbb{H}}^m$, one has 
\begin{displaymath}
{\widetilde{\tau }}_{-}\big( u +vj, z+w j \big) = 
(u^T +v^Tj) I_n (z + wj)^q. 
\end{displaymath}
This takes care of the case when $G = \Oo (\widetilde{V}, {\widetilde{\sigma }}_{-}, 
\widetilde{\tau })$. 

Suppose that $\widetilde{\tau }$ is \emph{alternating}. We may assume that the matrix of $\widetilde{\tau}$ with respect to the standard basis of ${\C}^n = {\C}^m \times {\C}^m$ is $J_n=${\tiny $\begin{pmatrix} 0 & I_m \\ -I_m & 0 \end{pmatrix}$}, that is, $\widetilde{\tau }\big( (z,w)^T, (z',w')^T \big) =$\raisebox{1pt}{ {\tiny  
 \raisebox{1pt}{${\begin{pmatrix} z' \\ w' \end{pmatrix}}^T$}}}$J_n $ \raisebox{2pt}{{\tiny $\! \! \begin{pmatrix} z\\ w \end{pmatrix}$}}. For $0 \le p \le m$ let  
\begin{displaymath}
{\widetilde{\sigma }}^{(p)}: {\C}^n = {\C}^m \times {\C}^m \rightarrow {\C}^n:
(z,w)^T \mapsto (-I_{m-p,p}\overline{w}, I_{m-p,p}\overline{z})^T. 
\end{displaymath}
Clearly ${\widetilde{\sigma }}^{(p)}_{-}$ is anti-linear and 
$({\widetilde{\sigma }}^{(p)}_{-})^2 = -{\mathrm{id}}_{{\C}^n}$. Replacing 
$\widetilde{\sigma }$ by $-\widetilde{\sigma}$, we may assume that $0 \le p \le m$. 
Let $\widetilde{u} = (z,w)^T$ and $\widetilde{v} = (z',w')^T$. Then 
\begin{align*}
\widetilde{\tau } \big( {\widetilde{\sigma }}^{(p)}_{-}(\widetilde{u}, 
{\widetilde{\sigma }}^{(p)}_{-}(\widetilde{v}) \big) & =   \\
&\hspace{-1in} = \left( \mbox{\footnotesize $\begin{pmatrix} 0 & -I_{m-p,p} \\ 
I_{m-p,p} & 0 \end{pmatrix} $} \widetilde {v} \right)^{\! \! T} 
\mbox{\footnotesize $\begin{pmatrix} 0 & I_m \\ - I_m & 0 \end{pmatrix} $} \left( \mbox{\footnotesize $\begin{pmatrix} 0 & -I_{m-p,p} \\ I_{m-p,p} & 0 \end{pmatrix} $} \widetilde{u} \right)  \\
& \hspace{-1in} = {\widetilde{v}}^{\, T} \mbox{\footnotesize  $\begin{pmatrix}
0 &  I_m \\ - I_m & 0  \end{pmatrix} $} \widetilde{u} 
 = \overline{\widetilde{\tau }\big( \widetilde{u}, \widetilde{v} \big)},  
\end{align*}
which implies
\begin{align*}
{\widetilde{\tau }}^{(p)}\big( \widetilde{u}, {\widetilde{\sigma }}^{(p)}_{-}(
\widetilde{v}) \big) & = 
- {\widetilde{\tau }}^{(p)}\big( {\widetilde{\sigma }}^{(p)}_{-}
({\widetilde{\sigma }}^{(p)}_{-}(\widetilde{u}) ), {\widetilde{\sigma }}^{(p)}_{-}(\widetilde{v}) \big) 
= -\overline{{\widetilde{\tau }}^{(p)}\big( {\widetilde{\sigma }}^{(p)}_{-}(\widetilde{u}, 
\widetilde{v}) \big)} . 
\end{align*}
So 
\begin{align*}
\big({{\widetilde{\tau }}^{(p)}_{-}( \widetilde{u}, \widetilde{v})} \big)^q & = 
\big( {\widetilde{\tau }}^{(p)}( \widetilde{u}, \widetilde{v}) + 
{\widetilde{\tau }}^{(p)} \big( \widetilde{u}, {\widetilde{\sigma }}^{(p)}_{-}(\widetilde{v}) \big) j \big)^q 
%\\ & 
= {\widetilde{\tau }}^{(p)} ( \widetilde{u}, \widetilde{v} ) 
- \overline{{\widetilde{\tau }}^{(p)} \big( \widetilde{u}, 
{\widetilde{\sigma }}^{(p)}_{-}(\widetilde{v}) \big) } j \\
& = {\widetilde{\tau }}^{(p)}\big( \widetilde{u}, \widetilde{v} \big) + 
{\widetilde{\tau }}^{(p)}\big( {\widetilde{\sigma }}^{(p)}_{-}(\widetilde{u}, 
\widetilde{v} )\big) j 
%\\ & 
= -{\widetilde{\tau }}^{(p)}( \widetilde{v}, \widetilde{u} ) - 
{\widetilde{\tau }}^{(p)}\big( \widetilde{v}, \widetilde{\sigma}^{(p)}_{-}(\widetilde{u}) \big) j \\
& = -{{\widetilde{\tau }}^{(p)}}_{-}(\widetilde{v}, \widetilde{u} ) ,
\end{align*}
that is, 
${\widetilde{\tau }}^{(p)}_{-}$ is \emph{hamiltonian alternating}. Let 
${\epsilon}' = {\{ \frac{\mathrm{i}}{\sqrt{2}}(e_{\ell }, f_{\ell }) \} }^m_{\ell =1}$ 
be a basis for ${\widetilde{V}}_{{\sigma }^{(p)}_{-}}$. Then 
\begin{align*}
{\widetilde{\tau }}^{(p)}_{-}\big( \ttfrac{\mathrm{i}}{\sqrt{2}}(e_{\ell },f_{\ell }), 
\ttfrac{\mathrm{i}}{\sqrt{2}}(e_k,f_k) \big) & =  \\
&\hspace{-1.25in} = -\onehalf {\widetilde{\tau }}^{(p)}\big( (e_{\ell }, f_{\ell }), (e_k,f_k) \big)  -\onehalf {\widetilde{\tau }}^{(p)}\big( (e_{\ell }, f_{\ell }), 
{\widetilde{\sigma }}^{(p)}_{-}(e_k,f_k) \big) j
\\
&\hspace{-1.25in}= -\onehalf {\widetilde{\tau }}^{(p)}\big( (e_{\ell }, 
f_{\ell }), (e_k,f_k) \big) -\onehalf {\widetilde{\tau }}^{(p)} \big( (e_{\ell }, f_{\ell }), 
{\delta }^{m-p,p}_k(-f_k, e_k) \big) j , 
\\
&\hspace{-.25in} \mbox{where ${\delta }^{m-p,p}_k =${\tiny $\left\{ \begin{array}{rl}
\hspace{-5pt} -1, & \hspace{-5pt} \mbox{if $1 \le k \le m-p$} \\ 
\hspace{-5pt} 1, & \hspace{-5pt} \mbox{if $m-p+1 \le k \le n$} \end{array} 
\right. $ }} \\
& \hspace{-1.25in} = {\delta }^{m-p,p}_k {\delta }_{\ell, k} \, j . 
\end{align*}
In other words, with respect to the basis ${\epsilon}'$ of 
${\widetilde{V}}_{{\sigma }^{(p)}_{-}}$ one has 
\begin{displaymath}
{\widetilde{\tau }}^{(p)}_{-}\big( u+vj, z+wj \big) = 
(u^T+v^Tj) (j\, I_{m-p,p}) (z+wj)^q .
\end{displaymath} 
This takes care of the case $\Spp (\widetilde{V}, {\widetilde{\sigma }}_{-}, 
{\widetilde{\tau}}^{(p)})$. \medskip 

\noindent \textsc{Case} $\mathbf{2.}$ \medskip 

Let $\widetilde{\sigma }: \widetilde{V} \rightarrow {\widetilde{V}}^{\ast }$ be 
an anti-linear mapping. Two such mappings $\widetilde{\sigma }$ and 
${\widetilde{\sigma }}'$ are equivalent if and only if ${\widetilde{\sigma }}' = 
\alpha k^{\ast} \widetilde{\sigma }k$ for some $\alpha \in \C \setminus \{ 0 \} $ 
and some $k \in \Gl (\widetilde{V})$. Here $k^{\ast }: {\widetilde{V}}^{\ast } 
\rightarrow {\widetilde{V}}^{\ast} $ is defined by 
$k^{\ast }({\widetilde{v}}^{\ast })(\widetilde{w}) = 
{\widetilde{v}}^{\ast }(k(\widetilde{w}))$ for every ${\widetilde{v}}^{\ast } \in 
{\widetilde{V}}^{\ast }$ and every $\widetilde{w} \in \widetilde{V}$. Let $\widetilde{V} = {\C}^n$. Set 
${\widetilde{\sigma }}^{(p)} : {\C }^n \rightarrow 
({\C }^n)^{\ast }: z \mapsto (I_{n-p,p}\overline{z})^T$. 
Then ${\widetilde{\sigma }}^{(p)}$ is an anti-linear mapping, since for every $z \in 
{\C}^n$ and every $\alpha \in \C$ one has 
${\widetilde{\sigma }}^{(p)}(\alpha z) = 
(\overline{\alpha z})^TI_{n-p,p}= \overline{\alpha }({\overline{z}}^TI_{n-p,p}) 
= \overline{\alpha } \, {\widetilde{\sigma }}^{(p)}(z)$. 
Let ${\widetilde{\tau }}^{(p)}_{\ast }: {\C }^n \times {\C }^n \rightarrow \C: 
(z,w) \mapsto \big( {\widetilde{\sigma }}^{(p)}(z)\big) (w) = 
{\overline{z}}^TI_{n-p,p}w$. Then ${\widetilde{\tau }}^{(p)}_{\ast }$ is a hermitian form, since 
\begin{align*}
\overline{{\widetilde{\tau }}^{(p)}_{\ast }(z,w)} & = 
\overline{ {\widetilde{\sigma }}^{(p)}(w)(z)} = \overline{{\overline{w}}^TI_{n-p,p}z} 
= w^TI_{n-p,p}\overline{z} \\
&  = {\overline{z}}^TI_{n-p,p}w  = {\widetilde{\sigma }}^{(p)}(w)(z) = 
{\widetilde{\tau }}^{(p)}_{\ast }(w,z). 
\end{align*}
This takes care of the case $G = \Gl (\widetilde{V}, {\widetilde{\tau }}^{(p)}_{\ast })$. 
\bigskip  

\noindent {\Large \textbf{Appendix 2. Indecomposable nilpotent types}} \bigskip 
%%%%%%%%%%%%%%%%%%%

In this appendix we give a table which lists all the indecomposable nilpotent 
types for the real nonexceptional Lie algebras. \medskip 

\noindent \hspace{.35in}\begin{tabular}{lll}
\multicolumn{1}{c}{Lie algebra} & \multicolumn{1}{c}{indecomposable type} & 
\multicolumn{1}{c}{index} \\ \hline
$\gl (V, {\sigma }_{+})$ & ${\Delta }_h(0)$ &  \\
$\gl (V, {\sigma }_{-})$ & ${\Delta }_h(0,0)$ &  \\
$\gl (V, {\tau }_{\ast })$ & ${\Delta }^{\varepsilon}_h(0)$, \mbox{$h$ even} &  
$\onehalf (h+1 - \delta )$ \\
\rule{0pt}{12pt}& ${\Delta }^{\varepsilon}_h(0)$, \mbox{$h$ odd} &  
$\onehalf (h+1)$ \\
\rule{0pt}{12pt}$\oo (V, {\sigma }_{+}, \tau )$  & ${\Delta }^{\varepsilon }_h(0)$, \mbox{$h$ even} & $\onehalf (h+1 - \delta  )$ \\
\rule{0pt}{12pt} \hspace{-2pt} & 
${\Delta }^{\varepsilon }_h(0,0)$, \mbox{$h$ odd} & $\onehalf (h+1)$ \\
\rule{0pt}{12pt}$\oo (V, {\sigma }_{-}, \tau )$ & ${\Delta }_h(0,0)$, \mbox{$h$ even}  & \\
\rule{0pt}{12pt} & ${\Delta }^{\varepsilon }_h(0,0)$, \mbox{$h$ odd}  
& \\
$\spp (V, {\sigma }_{+}, \tau )$ & ${\Delta }_h(0,0)$, \mbox{$h$ even} & \\
 & ${\Delta }^{\varepsilon}_h(0)$, \mbox{$h$ odd} & \\
$\spp (V, {\sigma }_{-}, \tau )$ & ${\Delta }^{\varepsilon }_h(0,0)$, 
\mbox{$h$ even} & $\onehalf (2(h+1) - \delta )$ \\
\rule{0pt}{12pt}& ${\Delta }^{\varepsilon }_h(0,0)$, \mbox{$h$ odd} & $h+1$ \\
\end{tabular} \bigskip

\noindent \hspace{.25in}\parbox[t]{4in}{Table A2. List of indecomposable nilpotent types. Here $\delta = (-1)^{h/2}\varepsilon $ and ${\varepsilon }^2 =1$.} \bigskip

Table A2 is taken from Table II p.349 of \cite{burgoyne-cushman}. \bigskip 
 
\noindent {\Large \textbf{Appendix 3. Indecomposable semisimple types }} \bigskip 
%%%%%%%%%%%%%%%%%%%

In this appendix we give a table which lists all the indecomposable semisimple 
types with $S=0$ for the real affine nonexceptional Lie algebras, see 
\cite[Table A p.361]{burgoyne-cushman}.\bigskip 

\noindent \hspace{.35in}\begin{tabular}{lcc} 
\multicolumn{1}{c}{Lie algebra} & \multicolumn{1}{c}{indecomposable type} & 
\multicolumn{1}{c}{comments} \\ \hline
$\gl (V)$ & ${\Delta }_0(0)$ &  \\
$\gl (V, {\sigma }_{+})$ & ${\Delta }_0(0)$ & (c) \\
$\gl (V, {\sigma }_{-})$ & ${\Delta }_0(0,0)$ & (b) \\
$\gl (V, {\tau}_{\ast})$ & ${\Delta }^{\varepsilon}_0(0)$ & \\
$\oo (V, \tau )$ & ${\Delta }_0(0)$ & \\ 
$\oo (V, {\sigma }_{+}, \tau )$ & ${\Delta }^{\varepsilon}_0(0)$ & 
$({\Delta }^{\varepsilon}_0(0))^{\mathrm{c}} = {\Delta }_0(0)$ \\
$\oo (V, {\sigma }_{-}, \tau )$ & ${\Delta }^{\varepsilon}_0(0)$ & (b) \\
$\spp (V, \tau )$ & ${\Delta }_0(0,0)$ & \\ 
$\spp (V, {\sigma }_{+}, \tau )$ & ${\Delta }_0(0,0)$ & (c) \\
$\spp (V, {\sigma }_{-}, \tau )$ & ${\Delta }^{\varepsilon}_0(0,0)$ & 
$({\Delta }^{\varepsilon}_0(0,0))^{\mathrm{c}} = {\Delta }_0(0,0)$
\end{tabular}\bigskip

\noindent \hspace{.25in}\parbox[t]{4.25in}{Table 3. List of indecomposable 
semisimple types with $S=0$.} \medskip 

Let $\Delta $ be an indecomposable semisimple type for $g(W, \sigma , \tau )$, 
represented by the pair $(W,S; \sigma , \tau )$. Then 
$(W,S; \tau ) \in {\Delta }^{\mathrm{c}}$ for $S \in g(W; \tau )$. Suppose that 
${\Delta }^{\mathrm{c}}_1$ is an indecomposable summand of 
${\Delta }^{\mathrm{c}}$, which is represented by the pair $(W_1, S; {\tau })$. 
Since ${\sigma }^{\ast}\tau = \overline{\tau }$, the pair $(\sigma (W_1), S; 
{\sigma }^{\ast }\tau )$ is well defined. Denote its type by 
$\sigma ({\Delta}^{\mathrm{c}}_1)$. Note that either $\sigma (W_1) = W_1$ or 
$W_1 \cap \sigma (W_1) = \{ 0 \} $. Since ${\sigma }^2= \pm 1$, $S=0$, and 
$\Delta $ is indecomposable, there are two possible cases: \smallskip 

(b) ${\Delta }^{\mathrm{c}} = 
{\Delta }^{\mathrm{c}}_1 + \sigma ({\Delta }^{\mathrm{c}}_1)$ 
and ${\Delta }^{\mathrm{c}}_1 = \sigma ({\Delta }^{\mathrm{c}}_1)$;   

(c) ${\Delta }^{\mathrm{c}} = {\Delta }^{\mathrm{c}}_1$ and ${\Delta }^{\mathrm{c}}_1 = \sigma ({\Delta }^{\mathrm{c}}_1)$.

\end{document}